\documentclass{amsart}
\usepackage{amssymb, epsfig}
\newtheorem{thm}{Theorem}[section]
\newtheorem{cor}[thm]{Corollary}
\newtheorem{lem}[thm]{Lemma}
\newtheorem{prop}[thm]{Proposition}

\newenvironment{acknowledgements}
{\subsubsection*{Acknowledgements}}{}
\renewcommand{\bibname}{\sc}
\renewcommand{\and}{\textnormal{and}\ }
\newenvironment{proof*}{\begin{proof}}{\end{proof}}

\newcommand{\scrg}{\mathcal{G}}
\newcommand{\scrc}{\mathcal{C}}
\newcommand{\scri}{\mathcal{I}}
\newcommand{\scrj}{\mathcal{J}}
\newcommand{\scrk}{\mathcal{K}}
\newcommand{\scrm}{\mathcal{M}}
\newcommand{\scrl}{\mathcal{L}}
\newcommand{\bbz}{\mathbb{Z}}
\newcommand{\bbr}{\mathbb{R}}
\newcommand{\bbq}{\mathbb{Q}}
\newcommand{\bbc}{\mathbb{C}}

\DeclareMathOperator{\id}{id}
\newcommand{\ttfrac}[2]{ {{#1}/{#2}} }
\newcommand{\tnfrac}[2]{ {({#1})/{#2}} }
\newcommand{\tdfrac}[2]{ {{#1}/({#2})} }
\newcommand{\ndfrac}[2]{ {({#1})/({#2})} }

\title{Polynomial splittings of Casson--Gordon invariants}
\author{Se-Goo Kim}
\address{Department of Mathematics,
University of California,
Santa Barbara, CA 93106, U.S.A.}
\email{sekim@math.ucsb.edu}
\urladdr{www.math.ucsb.edu/\~{}sekim}
\begin{document}

\begin{abstract}
In this paper we prove that the
Casson--Gordon invariants of the connected sum
of two knots split when the Alexander polynomials of the knots are
coprime.  As one application, for any knot $K$, all but finitely many
algebraically slice twisted doubles of $K$
are linearly independent in the knot concordance group.
\end{abstract}

\maketitle

%%%%%%%%%%%%%%%%%% Section %%%%%%%%%%%%%%%%%%%%%%%%%

\section{Introduction}
In his classification of the knot concordance groups,
Levine~\cite{lev69b} defined the algebraic concordance group,
$\scrg$, of Witt classes of Seifert matrices and a homomorphism
from the knot concordance group, $\scrc$, of knots in the
$3$-sphere $S^3$ to $\scrg$. Casson and Gordon~\cite{cg86} proved
that the kernel of Levine's homomorphism $\scrc\to\scrg$, the
concordance group of {algebraically slice knots}, is
nontrivial. Gilmer~\cite{gil93} used the work of~\cite{cg86} to
define a Witt type group $\Gamma^+$ and showed that there are
homomorphisms $\scrc\to\Gamma^+\to\scrg$.  The group $\Gamma^+$ is
roughly characterized by the property that a class of knots maps
to zero in $\Gamma^+$ if and only if all of Levine's invariants
and the Casson--Gordon invariants of a representative of the
class vanish. However, the definition of $\Gamma^+$ used here is
modified from the one used by Gilmer to correct for an error
in~\cite{gil93}, as described in Section~\ref{sec:def}.

It follows from Levine's work~\cite{lev69a} that
if the connected sum of two knots with relatively prime
Alexander polynomials maps to zero in $\scrg$, then so does each knot.
We show a similar result for the Casson--Gordon invariants as follows.
\begin{thm}
\label{thm:alexander}
Let $K_1$ and $K_2$ be knots with relatively prime Alexander polynomials
in $\bbq[t^{\pm1}]$.
Suppose that either $K_1$ or $K_2$ has a non-singular Seifert form.
Then if $K_1\# K_2$ is zero in $\Gamma^+$ then so are both $K_1$
and $K_2$.
\end{thm}

To demonstrate the strength of this result we study the family of
$k$-twisted doubles of a given knot
$K$, denoted $D_k(K)$.  This family contains an infinite number of
algebraically slice knots and these algebraically slice knots have
been the subject of careful study.  Casson and Gordon~\cite{cg86} found
the first
examples of nontrivial concordance classes in the kernel of Levine's
homomorphism using the family $D_k(U)$ where $U$ is the unknot.
Since then \cite{jia81,lit84,ln99a, ln99b,tam99} have
found infinite linearly independent families of algebraically slice knots
among the knots $D_k(U)$.
In each case these families were very scarce: roughly one knot was
chosen
for each prime integer.  Theorem~\ref{thm:alexander} will   yield that
for every knot $K$ (not just the unknot $U$)  the set of all
algebraically slice knots in the family of knots $D_k(K)$ is (with finite
exceptions) linearly independent.  More precisely:

\begin{thm}
\label{thm:independent}
(a) For any knot $K$, all but finitely
many algebraically slice twisted doubles of $K$ are linearly
independent in~$\Gamma^+$ and so in~$\scrc$.

(b) If $\sigma_r(K)\ge 0$ for all $r$, then all algebraically slice
twisted doubles of $K$ except the untwisted and 2-twisted ones
are linearly independent in~$\Gamma^+$ and so in~$\scrc$, where
$\sigma_r$ denotes the averaged Tristram--Levine signature.  
In addition, if $\sigma_r(K)>0$ for some $r$, then
all algebraically slice
twisted doubles of $K$ except the untwisted one
are linearly independent in~$\Gamma^+$ and so in~$\scrc$.
\end{thm}

\begin{cor}
\label{cor:independent}
(a) All algebraically slice twisted doubles of the
unknot except the two known to be slice
(untwisted and 2-twisted ones) are linearly
independent in the knot concordance group.

(b) There are infinitely many knots $K$ for each of which all 
algebraically slice twisted doubles of $K$ except
the untwisted one are linearly independent in
the knot concordance group.
\end{cor}

The full set of knots $D_k(K)$ contains knots which are not
algebraically slice, representing elements of infinite order, order 2,
and order 4 in  $\scrg$.  It has not been possible to prove that this
set is linearly independent in $\scrc$, but we do have the following
theorem, based on recent work of Friedl~\cite{fri03}:

\begin{thm}
\label{thm:ribbon}
(a) For any knot $K$, there is a set
$\scrk$ containing all twisted doubles of $K$ except
a finite number of nonnegatively twisted ones such that
no nontrivial linear combinations of elements in $\scrk$ are ribbon.

(b) If $\sigma_r(K)\ge 0$ for all $r$, then
no nontrivial linear combinations of twisted doubles of $K$ except those
with $0,1,2$ twists are ribbon.
In addition, if $\sigma_r(K) > 0$
for $r=\frac{2}{9},\frac{1}{3},\frac{2}{5}$, then
no nontrivial linear combinations of twisted doubles of $K$ except the
untwisted one are ribbon.
\end{thm}

In the past,
the construction of independent knots
depended on finding knots for which some branched covers had homology
groups of order divisible by distinct primes.  Such an approach could
conceivably work with doubled knots by using high degree covers, but the
argument would be far more burdensome than the one we give.
A paper in preparation will address another application of
Theorem~\ref{thm:alexander} that there are
examples of linearly independent algebraically slice
knots having the same homology on all prime
power fold branched covers, in which case
no such approach could possibly work, and
our approach using the splitting associated with the polynomial is
definitely required.

The paper is organized as follows: 
In Section~\ref{sec:def} we summarize the modified results
of~\cite{gil93} about the Casson--Gordon invariants on slice knots and
ribbon knots to correct for an error in~\cite{gil93}. In
Section~\ref{sec:alexander} we prove Theorem~\ref{thm:alexander}. In
Section~\ref{sec:double} we state that all algebraically slice twisted
doubles of a knot (with finite exceptions) have infinite order in
$\Gamma^+$ and similar results concerning their ribbonness.  We also
summarize some facts on the Casson--Gordon invariants of genus 1 knots
and the Tristram--Levine signatures of satellite knots and torus knots.
In Section~\ref{sec:a-slice} we estimate the Casson--Gordon invariants of
algebraically slice $D_k(K)$ and prove Theorem~\ref{thm:independent} and
Corollary~\ref{cor:independent}. In Section~\ref{sec:ribbon} we estimate
the Casson--Gordon invariants of all twisted doubles $D_k(K)$ for double
branched covers and prove Theorem~\ref{thm:ribbon}.

%%%%%%%%%%%%%%%%%% Section %%%%%%%%%%%%%%%%%%%%%%%%%

\section{Gilmer's obstructions}
\label{sec:def}
In this section we state the modified results of~\cite{gil93} about
the Casson--Gordon invariants on slice knots and ribbon knots to correct
for an error in~\cite{gil93}. Conventions of~\cite{nai96} are followed
here rather than those of~\cite{gil93}
since we will use the formula for the Casson--Gordon invariants for genus
1 knots given in~\cite{nai96}. 

Let $K$ be a knot in the $3$-sphere $S^3$ with Seifert surface $F$
having intersection pairing $\langle\;,\,\rangle$.
Let $i_+\!:H_1(F)\to H_1(S^3-F)$
denote the map which pushes a class off in the positive normal direction
and let $i_-$ denote the map given by pushing off the other way.

Let $\theta$ denote the Seifert pairing, i.e.~
$\theta(x,y)= lk(x,i_+ y)$.
Let $A$ denote the Seifert matrix for $\theta$ with respect
to some basis $\{a_1,\ldots,a_{2g}\}$ for $H_1(F)$.  With respect to this
basis the intersection form on $F$ is then given by the matrix $A^t-A$.
Let $\{\alpha_1,\ldots,\alpha_{2g}\}$ denote the basis for $H_1(S^3-F)$ such
that $lk(a_i,\alpha_j)=\delta^i_j$ as in~\cite[page~209]{rol76}.
Then $i_+$ with respect to these bases is given by $A$.  Define
$j\!:H_1(S^3-F)\to H_1(F)$ by $\langle jx,y\rangle=lk(x,y)$, so $j$ is given
by the matrix $(A-A^t)^{-1}$ with respect to the above bases.

Following~\cite{ker71,sto77},
we define \emph{the associated isometric
structure} $s\!:H_1(F)\to H_1(F)$ by the equation
$\theta(x,y)=\langle sx,y\rangle$.  We see that $s$ is $j\circ i_+$ and
is given by the matrix $G=(A-A^t)^{-1}A$.  Note $s-1$ is given by
$G-I=(A-A^t)^{-1}A-I=(A-A^t)^{-1}A^t$ and is actually
$j\circ i_-$.
Let $M^q$
denote the $q$-fold branched cyclic cover of $S^3$ along $K$.
Then Seifert showed that
$G^q-(G-I)^q$ is a presentation matrix for $H_1(M^q)$ 
(for a more recent reference, see \cite[lemma~1]{gil93}).

We are interested in $H^1(M^q;\bbq/\bbz)$, the set of characters on
$H_1(M^q)$. Define $\varepsilon^q$ to be the endomorphism of
$H_1(F)$ given by $s^q-(s-1)^q$ and $N^q\subset H_1(F;\bbq/\bbz)$
to be the kernel of $\varepsilon^q\otimes \id_{\bbq/\bbz}$.
Gilmer~\cite{gil93} proved that $H^1(M^q;\bbq/\bbz)$ is isomorphic
to $N^q$ and the isomorphism can be uniquely constructed up to
covering translations.  So we may view the Casson--Gordon
invariants $\tau(K,\,)$ as a function on $N^q$. See~\cite{cg86}
for the definition of the Casson--Gordon invariants.  From now on
$q$ will always denote a power of a prime.  For a prime
$p$, let $N^q_p$ denote the $p$-primary component of $N^q$.

A Seifert form $\theta$ on $H_1(F)$ is said to be \emph{null-concordant}
if there is an $s$-invariant direct summand $Z$ of $H_1(F)$ such that
$Z=Z^\bot$ with respect to the intersection pairing
$\langle\;,\,\rangle$.
Such a direct summand $Z$ is called a \emph{metabolizer} for the isometric
structure $s$ associated to $\theta$.

We say a knot $K$ is \emph{slice} if $K$ bounds a smoothly
embedded 2-disk $D$ in the $4$-ball $B^4$ with
$\partial(B^4,D)=(S^3,K)$. Knots $K_1$ and $K_2$ are called
\emph{concordant} if $K_1\# -K_2$ is slice, where $-K$ denotes
the mirror image of $K$ with reversed orientation.  The set of
concordance classes of knots forms an abelian group under
connected sum, called the \emph{knot concordance group}
and denoted~$\scrc$.
We remark here that the result all apply in the category of locally flat
oriented manifolds and pairs as well.

We say a knot $K$ is \emph{ribbon} if $K$ bounds a smoothly immersed
2-disc $f(D^2)$ in $S^3$ for a smooth map $f\!:D^2\to S^3$
which has the property that each component of self-intersection is an arc
$A\subset f(D^2)$ for which $f^{-1}(A)$ is two arcs in~$D^2$, one of
which lies in the interior of $D^2$.
It is easy to see that all ribbon knots are slice.
There is no known concept of the knot \emph{ribbon} concordance group
more along the line of the knot concordance group as defined above.  The
difficulty is that it is unknown whether the following is true:
If knots $K$ and $K\# J$ are ribbon, then $J$ is ribbon.  In fact, this
is an equivalent statement of Fox's conjecture: All slice knots are
ribbon.

Gilmer~\cite{gil83, gil93} combined the slicing obstructions
of~\cite{lev69b} with those of~\cite{cg86} in a
nontrivial way.  Recently Gilmer has announced that there is an error
in~\cite{gil83, gil93} (cf. \cite{fri03}).  However, the following two
weaker statements are known to be valid. Note that the first
statement has a weaker conclusion: the phrase ``all primes $p$''
is replaced by ``all but finitely many primes $p$''; the second
statement has a stronger hypothesis: ``a slice knot $K$'' is
replaced by ``a ribbon knot $K$.''
The first statement directly follows from Gilmer's original proof if
primes $p$ are chosen so that $p$ do not divide
the order of torsion of $H_1(R)$ (see~\cite{gil93} for the definition
of~$R$). The second statement is a corollary of Friedl's recent
work~\cite[theorem 8.3 and corollary 8.5]{fri03}.

\begin{thm}
\label{thm:gilmer}
(a) If $F$ is a Seifert surface
for a slice knot $K$ then there is a metabolizer $Z$ for the
isometric structure on $H_1(F)$ such that
$\tau(K,N^q_p\cap(Z\otimes\bbq/\bbz))$ vanishes for all prime
powers $q$ and all but finitely many primes $p$.

(b) If $F$ is a Seifert surface for a ribbon knot $K$
then there is a metabolizer $Z$ for the isometric structure on
$H_1(F)$ such that $\tau(K,N^q_p\cap(Z\otimes\bbq/\bbz))$ vanishes
for all prime powers $q$ and all primes $p$.
\end{thm}

We say that a knot $K$ has
\emph{vanishing Gilmer slice (resp.~ribbon) obstruction}
if the conclusion of Theorem~\ref{thm:gilmer}(a) (resp.~(b))
is satisfied for any Seifert surface $F$ of $K$.
Otherwise, we say that $K$ has \emph{nonvanishing}
Gilmer slice (resp.~ribbon) obstruction.

Gilmer~\cite{gil93} defined a Witt type group $\Gamma^+$ and a
homomorphism from the knot concordance group $\scrc$ to $\Gamma^+$
such that the class of a knot maps to zero if and only if it
satisfies the conclusion of Theorem~\ref{thm:gilmer}(b).
It became unknown whether the original $\Gamma^+$ is a group
and whether there is a homomorphism $\scrc\to\Gamma^+$
since the proof of the cancellation lemma of \cite[lemma 5]{gil93}
had a similar gap. 

On the other hand,
if the definition of $\Gamma^+$ is modified so that 
the class of a knot maps to zero if and only if it
satisfies the conclusion of Theorem~\ref{thm:gilmer}(a), then
$\Gamma^+$ is a group and there are still homomorphisms
$\scrc \to\Gamma^+\to\scrg$.
For instance,
using the notation in~\cite[page 12]{gil93}, the definition
of a metabolizer for $(U, \langle\,,\,\rangle,s,\tau^q_p)$ should be
modified as done in the conclusion of Theorem~\ref{thm:gilmer}(a), namely,
``all primes $p$'' should be replaced by ``all but finitely many
primes $p$.''   With this new definition of
metabolizer we can eliminate the gap in the proof of 
\cite[lemma 5]{gil93} and hence this new $\Gamma^+$ becomes
a group. Throughout this paper, $\Gamma^+$ is this modified one.

%%%%%%%%%%%%%%%%%% Section %%%%%%%%%%%%%%%%%%%%%%%%%

\section{Alexander polynomial and proof of Theorem~\ref{thm:alexander}}
\label{sec:alexander}
Let $K$ be a knot in the $3$-sphere $S^3$ with Seifert surface $F$,
Seifert pairing $\theta$, and Seifert matrix $A$.
We define $\Delta_K(t)=\det(A-tA^t)$, called the \emph{Alexander
polynomial of $K$ corresponding to the Seifert surface $F$}.
As is well known, the Alexander polynomial of a knot is
uniquely determined up to multiplication by $\pm t^n$ in $\bbq[t^{\pm1}]$.
Observe that
the characteristic polynomial for the associated isometric structure
$s$ to $K$ is $\det(xI-G)=\pm x^{2g}\det(A-(1-x^{-1})A^t)=
\pm x^{2g}\Delta_K(1-x^{-1})$, where $g$ is the genus of $F$.

To prove Theorem~\ref{thm:alexander}
we need a generalized notion of Seifert form.
Consider integral valued bilinear forms $\theta$ on finitely generated
free $\bbz$-modules $H$.  Define the \emph{transpose} of $\theta$,
denoted $\theta^t$, by
$\theta^t(x,y)=\theta(y,x)$ for all $x$ and $y$ in $H$.
We say that $\theta$ is
a \emph{Seifert form} if the form $\theta-\theta^t$ is unimodular,
i.e.~ the associated map $H\to \hbox{Hom}(H,\bbz)$, defined by
$x\mapsto (\theta-\theta^t)(x,\,)$, is an isomorphism.  A form is called
\emph{non-singular} if its associated map is injective.

The notions of isometric structure, metabolizer, and Alexander polynomial
extend to (algebraic) Seifert forms.  Observe that if $\theta$ is a
Seifert form on $H$, then the rank of $H$ must be even.
Polynomials in $\bbq[t^{\pm1}]$ are said to be
\emph{relatively prime} if their greatest common divisor is a unit.

The following lemma is a refinement of~\cite[proposition 3]{ker71}.  Its
geometrical origins are in the work of~\cite{lev69a} on the knot
concordance group.  They did not need the splitting of a metabolizer as
stated below, while we will need the splitting later.

\begin{lem}
\label{lem:alexander}
Let $\theta_1$ and $\theta_2$ be Seifert forms on $H_1$ and
$H_2$.  Suppose that their Alexander polynomials are relatively prime in
$\bbq[t^{\pm1}]$ and that either $\theta_1$ or $\theta_2$ is non-singular.
Then if $\theta_1\oplus\theta_2$ is null-concordant with a
metabolizer $Z$ for the associated isometric structure, then $\theta_1$ and
$\theta_2$ are null-concordant with metabolizers $Z_1$ and $Z_2$ for the
associated isometric structures, respectively, such that
$Z_i=Z\cap H_i,\ i=1,2$, and $Z=Z_1\oplus Z_2$.
\end{lem}

\begin{proof*}
Consider the associated isometric structures $s_i$ to $\theta_i$, $i=1,2$.
Then $s_1\oplus s_2$ on $H=H_1\oplus H_2$ is
the associated isometric structure $s$ to $\theta_1\oplus\theta_2.$
Let $Z_i=Z\cap H_i$, $i=1,2$, and let
$\varphi_i(x)=x^{2g_i}\Delta_{\theta_i}(1-x^{-1})$, where $2g_i$ is the
rank of $H_i$.
Since $\Delta_{\theta_1}$ and $\Delta_{\theta_2}$ are relatively prime in
$\bbq[t^{\pm1}]$ and since either $\theta_1$ or $\theta_2$ is non-singular,
$\Delta_{\theta_1}$ and $\Delta_{\theta_2}$ are relatively prime in
$\bbq[t]$.  Then $\varphi_1$ and $\varphi_2$ are also relatively prime in
$\bbq[x]$.  For, if
$f(x)$ is a common factor of $\varphi_1$ and $\varphi_2$,  then
$(1-t)^df\left({1}/{(1-t)}\right)$, $d=\deg f$, is a common factor
of
$\Delta_{\theta_1}$ and $\Delta_{\theta_2}$.
Thus there are polynomials $u_1$ and $u_2$ in $\bbz[x]$ and a non-zero
integer $c$ such that $u_1\varphi_1+u_2\varphi_2=c$.

For $z\in Z$, there are $z_1\in H_1$ and $z_2\in H_2$ with $z=z_1+z_2$.
As stated right after the definition of the Alexander polynomial, each
$\varphi_i$ is the characteristic polynomial for $s_i$, and hence
$\varphi_i(s_i)=0$.
Using this and $s(z_i)=s_i(z_i)$, we have
\begin{eqnarray*}
cz_1 & =  & u_1(s)\varphi_1(s)z_1+u_2(s)\varphi_2(s)z_1 \\
& = & u_1(s_1)\varphi_1(s_1)z_1+u_2(s_1)\varphi_2(s_1)z_1 \\
& = & u_2(s_1)\varphi_2(s_1)z_1 \\
& = & u_2(s_1)\varphi_2(s_1)z_1+u_2(s_2)\varphi_2(s_2)z_2 \\
& = & u_2(s)\varphi_2(s)z.
\end{eqnarray*}
Since $Z$ is $s$-invariant, $cz_1=u_2(s)\varphi_2(s)z\in Z$.
Since $Z$ is a direct summand of $H$,
this implies $z_1 \in Z$, and hence $z_1\in Z_1$.
Similarly, $z_2\in Z_2$.  So $Z=Z_1+Z_2$.  This now implies that
$Z=Z_1\oplus Z_2$ since $H=H_1\oplus H_2$ and $Z_i=Z\cap H_i$.

Since $Z$ is $s$-invariant, each $Z_i$ is $s_i$-invariant.
Since $H/Z = H_1/Z_1\oplus H_2/Z_2$ is torsion free, each $Z_i$ is a
direct summand of $H_i$.
Since the intersection pairing $\langle\;,\,\rangle$ on $H$ is unimodular,
$Z_i=Z_i^\perp$ on $H_i$.
Thus each $Z_i$ is a metabolizer for $s_i$.
\end{proof*}

\begin{proof*}[Proof of Theorem~\ref{thm:alexander}.]
Let $F_1$ and $F_2$ be Seifert surfaces for $K_1$ and $K_2$.  Then a
boundary connected sum $F_1\natural F_2$ is a Seifert surface for
$K_1\#K_2$.  Let $Z$ be a metabolizer for the isometric structure on
$H_1(F_1\natural F_2)=H_1(F_1)\oplus H_1(F_2)$ satisfying the
conclusion of Theorem~\ref{thm:gilmer}(a) with the exceptional primes
$p_1,\ldots,p_n$, i.e.~
$\tau(K_1\# K_2, N^q_{p}\cap(Z\otimes\bbq/\bbz))$ vanishes
for all prime powers $q$ and all primes $p$ except $p_1,\ldots,p_n$.

Then by Lemma~\ref{lem:alexander}
there are metabolizers $Z_1$ and
$Z_2$ for the isometric structures on $H_1(F_1)$ and $H_1(F_2)$,
respectively, with $Z=Z_1\oplus Z_2$.
Let $q$ be a power of a prime.
Let $N=\ker\varepsilon^q\otimes \id_{\bbq/\bbz}$ and
$N_i=\ker\varepsilon^q_{i}\otimes \id_{\bbq/\bbz}$, where $\varepsilon^q$
and $\varepsilon^q_{i}$ are the endomorphisms of $H$ and $H_i$,
respectively, as denoted in Section~\ref{sec:def}.  Then since
$\varepsilon^q=\varepsilon^q_{1}\oplus\varepsilon^q_{2}$, $N=N_1\oplus
N_2$, and $Z=Z_1\oplus Z_2$,
\[
N\cap(Z\otimes\bbq/\bbz) = (N_1\cap Z_1\otimes\bbq/\bbz) \oplus
(N_2\cap Z_2\otimes\bbq/\bbz).
\]

Let $N_p$ and $N_{i,p}$ denote the $p$-primary components of $N$
and $N_i$, respectively.  Let $\chi_1\in
N_{1,p}\cap(Z_1\otimes\bbq/\bbz)$. Then $\chi=\chi_1\oplus 0$ is
an element in $N_p\cap (Z\otimes\bbq/\bbz)$, where $0$ stands for
the trivial character in $N_2\cap (Z_2\otimes\bbq/\bbz)$.  By the
additivity of Casson--Gordon invariants~\cite[page~335]{lit84},
for all primes $p$ except
$p_1,\ldots,p_n$,
\[
0=\tau(K_1\# K_2,\chi)=\tau(K_1,\chi_1)+\tau(K_2,0).
\]
Also,  by \cite[corollary~B2]{lit84}
$\tau$ is determined by the algebraic concordance class of the
knot if the character is trivial.   This implies that $\tau(K_2,0)
= 0$ and hence $\tau(K_1,\chi_1)=0$.  Since $\chi_1$ was
chosen arbitrarily, we just have found a metabolizer $Z_1$ for the
isometric structure on $H_1(F_1)$ such that
$\tau(K_1,N_{1,p}\cap(Z_1\otimes\bbq/\bbz))$ vanishes for all
prime powers $q$ and all primes $p$ except $p_1,\ldots,p_n$,
i.e.~ $K_1$ is zero in $\Gamma^+$.  Similarly, $K_2$ is
zero in $\Gamma^+$. This completes the proof.
\end{proof*}

The proof given above also works for the Gilmer ribbon obstructions, in
which case there are no exceptional primes $p_1,\dots,p_n$.

\begin{cor}
\label{cor:ralexander}
Under the same conditions as in Theorem~\ref{thm:alexander},
if $K_1\# K_2$ has vanishing Gilmer ribbon obstruction, then
so do both $K_1$ and $K_2$.
\end{cor}

%%%%%%%%%%%%%%%%%% Section %%%%%%%%%%%%%%%%%%%%%%%%%

\section{Twisted doubles of a knot}
\label{sec:double}
In this section we state that all
but finitely many algebraically slice twisted doubles of a knot have
infinite order in the knot concordance group $\scrc$, in fact, in
$\Gamma^+$. We also state similar results concerning ribbonness in the
line of Theorem~\ref{thm:gilmer}(b).  The proofs will be given in the next
two sections.  In preparations, we also
summarize some facts on the Casson--Gordon invariants of genus 1 knots
and the Tristram--Levine signatures of satellite knots and torus knots.

\begin{figure}
\vspace{1em}
\setlength{\unitlength}{0.4pt}
\begin{picture}(314,180)
\put(-4,-2.5){\includegraphics[scale=0.4]{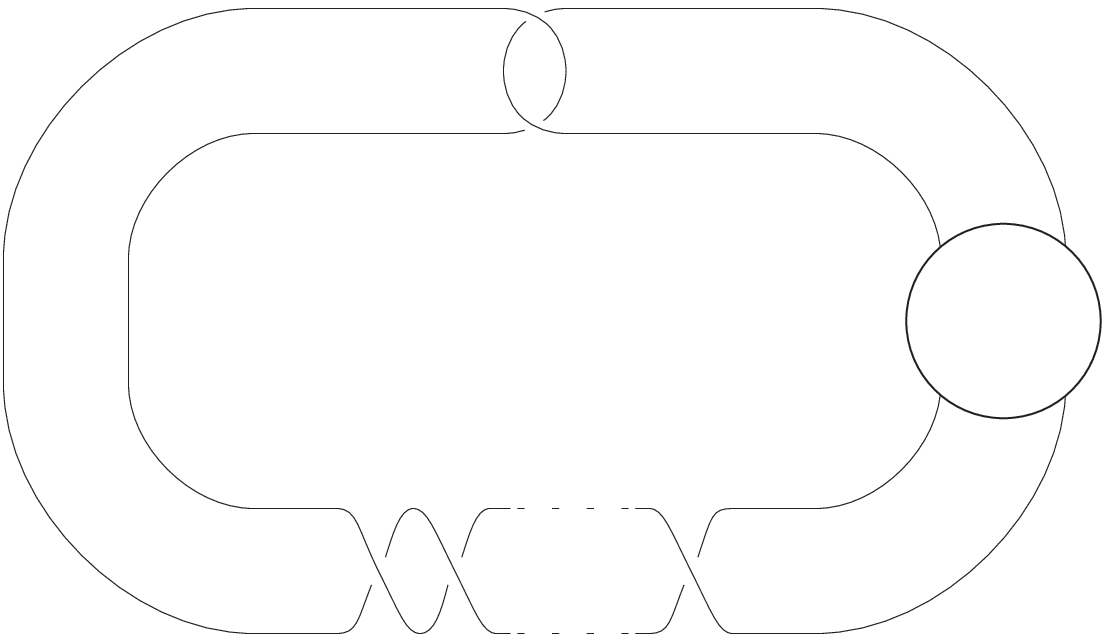}}
\put(154,45){\makebox(0,0)[b]{\tiny $k$ full twists}}
\put(288,90){\makebox(0,0){$K$}}
\end{picture}
\qquad
\setlength{\unitlength}{0.4pt}
\begin{picture}(464,180)
\put(-4,-2.5){\includegraphics[scale=0.4]{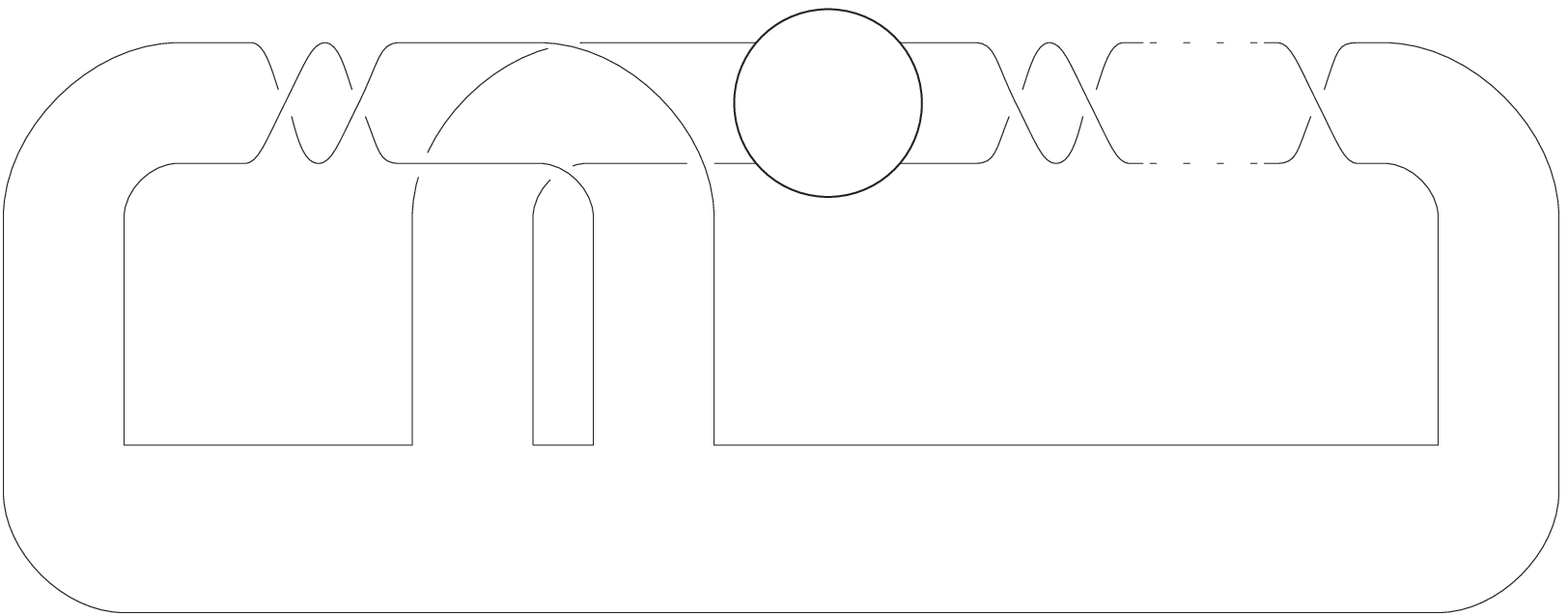}}
\put(342,122){\makebox(0,0)[t]{\tiny $k$ full twists}}
\put(246,152){\makebox(0,0){$K$}}
\end{picture}
\vspace{1em}
\caption{The $k$-twisted double of a knot $K$}
\label{fig:double}
\end{figure}

Let $K$ be a knot in the $3$-sphere $S^3$.
Let $D_k(K)$ denote the $k$-twisted double of $K$
as illustrated in Figure~\ref{fig:double}.
Here, $k$ may be negative.

%%%%%%%%%%%%%%%% Figure %%%%%%%%%%%%%%%%%%%%

The following theorem is due to~\cite[corollary~23]{lev69a}.

\begin{thm}
\label{thm:levine}
The $k$-twisted double of a knot $K$ is:

(a) of infinite order in the algebraic concordance group,
$\scrg$, if $k<0$;

(b) algebraically slice if $k\ge 0$ and $4k+1$ is a perfect square;

(c) of order $2$ in $\scrg$ if $k>0$, $4k+1$ is not a perfect
square, and every prime congruent to $3$ \textup{mod} $4$ has even
exponent in the prime power factorization of $4k+1$;

(d) of order $4$ in $\scrg$ if $k>0$ and some prime congruent to
$3$ \textup{mod} $4$ has odd exponent in the
prime power factorization of $4k+1$.
\end{thm}

Immediate corollaries are that $D_k(K)$ is
algebraically slice if and only if $k=l(l+1)$ for an integer $l\ge
0$ and that $D_k(K)$ has infinite order in
$\Gamma^+$ if $k<0$.

We will further prove the following two theorems.
For a nonnegative integer~$n$ let $nD_k(K)$ denote the connected sum
of $n$ copies of $D_k(K)$.
The first of the following two theorems
concerns the order of algebraically slice $D_k(K)$
in $\Gamma^+$.
The second concerns the
Gilmer ribbon obstructions of $nD_k(K)$ for not only algebraically
slice but all twisted doubles $D_k(K)$.

\begin{thm}
\label{thm:infinite}
(a) For any knot $K$, the algebraically slice
$k$-twisted double $D_{k}(K)$ has infinite order in $\Gamma^+$ for
all but finitely many $k$.

(b) If $\sigma_r(K)\ge 0$ for all $r$, the algebraically slice $D_k(K)$
has infinite order in $\Gamma^+$ for any $k\ne 0,2$, where $\sigma_r(K)$
denotes the averaged Tristram--Levine signature of $K$ (details will be
given later).  
In addition, if $\sigma_r(K)>0$ for some $r$, then
$D_2(K)$ has infinite order in $\Gamma^+$ as well.
\end{thm}

\begin{thm}
\label{thm:rinfinite} 
(a) For any knot $K$, there is a set
$\scri$ of all integers except a finite number of nonnegative integers
such that, for any $k\in\scri$ and any integer~$n\ne 0$,
$nD_k(K)$ has nonvanishing Gilmer ribbon obstruction.

(b) If $\sigma_r(K)\ge 0$ for all $r$, for any integer $n\ne 0$ and
any integer $k\ne 0,1,2$, $nD_k(K)$ has nonvanishing Gilmer ribbon
obstruction.  In addition, if $\sigma_r(K)
> 0$ for $r=\frac{2}{9}, \frac{1}{3}, \frac{2}{5}$, for any integer
$n\ne 0$, $nD_1(K)$ and $nD_2(K)$ have nonvanishing Gilmer ribbon
obstruction as well.
\end{thm}

We devote the remaining two sections to proving these theorems.
Before that, we summarize some useful facts in the rest of this section.

%%%%%%%%%%%%%%% Subsection %%%%%%%%%%%%%%%%%%

\subsection{Casson--Gordon invariants of a genus $1$ knot}
\label{sec:genus1}
We state the work of \cite[theorem~7]{nai96} that gives a formula for
$\tau$ for genus 1 knots in terms of the classical signatures.
We remark that \cite{gil83,gil93} first found the formula for the
2-fold branched cover case and algebraically slice case.

For a knot $K$ and a character $\chi$, $\tau(K,\chi)$ is defined to be an
element of the Witt group
$W(\bbc(t),\scrj)\otimes_\bbz \bbq$, where $\scrj$ denotes the involution
on $\bbc(t)$ given by complex conjugation and by the map $t\mapsto
t^{-1}$ and $W(\bbc(t),\scrj)$ is the Witt group of finite dimensional
hermitian inner product spaces.  For details, see~\cite{cg86}.
Let $W(\bbr)$ denote the Witt group of
finite dimensional inner product spaces over $\bbr$.  The signature
function $\sigma\!:W(\bbr)\to\bbz$ is an isomorphism.  Also there is a
natural map $W(\bbr)\to W(\bbc(t),\scrj)$ given by tensoring with $\bbc(t)$
over $\bbr$.  Composing this map with $\sigma^{-1}$ tensored with $\bbq$
gives a homomorphism $\rho\!:\bbq\to W(\bbc(t),\scrj)\otimes_\bbz
\bbq$.  Note that for each complex number $\zeta$ with $|\zeta|=1$, there
is a homomorphism $\sigma_\zeta\!:W(\bbc(t),\scrj)\otimes_\bbz
\bbq\to\bbq$ (see~\cite{cg86}).  It is easy to see that
$\sigma_1\circ\rho$ is the identity.

For any real number $r$, define
$A_r(K) = \left(1-e^{2\pi ir}\right)A
+ \left(1-e^{-2\pi ir}\right)A^t$, where $A$ is a Seifert
matrix of $K$, and define
$\sigma_r(K)$ to be $\sigma(A_r)$ if $A_r$ is non-singular and elsewhere
to be the average of the one-sided limits of $\sigma(A_r)$.  This
$\sigma_r(K)$ is a concordance invariant and
is equal to the Tristram--Levine signature~\cite{lev69b,tri69} of $K$
except perhaps at finitely many $r$.  It is an immediate consequence of
the definition that
$\sigma_r(K)=\sigma_{1-r}(K)$.
Thus we only need to consider $\sigma_r(K)$
for $0\le r\le \tfrac{1}{2}$.

The following is due to~\cite{nai96}. The case for $q=2$
and the case for algebraically slice knots are due to~\cite{gil83, gil93}.

\begin{thm} 
\label{thm:naik}
Let $F$ be a genus one Seifert surface of a knot $K$ and let
$A = \big(\begin{smallmatrix} a & -m \\ -(m+1) & b
\end{smallmatrix}\big)$ be the Seifert matrix of $K$ with respect
to a basis $\{x,y\}$ of $H_1(F)$.
Let $q$ and $d$ be powers of primes, $s$ an integer relatively prime to
$d$, and $N^q=\ker ((G^q-(G-I)^q)\otimes \id_{\bbq/\bbz})
\subset H_1(F;\bbq/\bbz)$, where
$G=(A-A^t)^{-1}A$. 
Suppose that $x\otimes {s}/{d}\in N^q$ and $d\mid a$.
Then the multiplicative inverse, $m^\ast$, of $m$ \textup{mod} $d$ exists
and
$x\otimes {s}/{d}$ defines a character $\chi$ for which
\[
\tau(K,\chi)=\rho \sum^{q-1}_{i=0}
\left(\sigma_{\frac{s_i}{d}}(J_x) + \frac{2(d-s_i)s_ia}{d^2}
-\sigma_{\frac{i}{q}}(K)\right),
\]
where $J_x$ is a simple closed curve on $F$ representing $x$ and,
for $i=0,\ldots, q-1$, $s_i$ is an integer such that $0<s_i<d$ and
$s_i\equiv (1+m^\ast)^i s$ \textup{mod} $d$.

In particular, 
(a) if $q=2$, then
\[
\tau(K,\chi)=
\rho\left(2\sigma_\frac{s_0}{d}(J_x)+\frac{4(d-s_0)s_0a}{d^2}
-\sigma_\frac{1}{2}(K)\right).
\]

(b) If $a=0$, then
$d\mid (m+1)^q-m^q$ implies $x\otimes {s}/{d}\in N^q$ and
\[
\tau(K,\chi)=\rho \sum^{q-1}_{i=0} \sigma_{\frac{s_i}{d}}(J_x).
\]
\end{thm}

%%%%%%%%%%%%%% Subsection %%%%%%%%%%%%%%%%%%%

\subsection{Satellite knots and torus knots}
Let $K$ be a knot in $S^3$.  By an \emph{axis} for $K$ of \emph{winding}
number
$w$ we mean an unknotted simple closed curve $\gamma$ in $S^3-K$ having
linking number $w$ with $K$.  Let $V$ be a solid torus complementary to a
tubular neighborhood of $\gamma$, with $K$ contained in the interior of
$V$.  There is a preferred generator $v$ for $H_1(V)$, specified by the
condition $lk(v,\gamma)=+1$.  For any knot $C$ in $S^3$ there is an
untwisted orientation-preserving embedding $h\!:V\to S^3$ taking $V$ onto
a tubular neighborhood of $C$ such that $C$ represents $h_\ast(v)$ in
$H_1(hV)$.  We say that the knot $h(K)$, denoted $C(K)$, is a
\emph{satellite} of $C$ with \emph{orbit} $K$, \emph{axis} $\gamma$, and
\emph{winding number} $w$.

The following is~\cite[theorem~2]{lit79}.

\begin{thm} 
\label{thm:satellite}
Let $C(K)$ be a satellite of $C$ with orbit $K$ and winding number~$w$.
Then
\[\sigma_r(C(K))=\sigma_{wr}(C)+\sigma_r(K).\]
\end{thm}

Let $T_{m,n}$ denote the $(m,n)$ torus link.
To fix orientation conventions $T_{2,2}$ is the positive Hopf link.
We will use the work of~\cite{lit79} on the signatures of
$T_{m,n}$ to prove the following proposition.
For a real number $z$,
let $[z]$ denote the greatest integer that is less than or equal to $z$.

\begin{prop}
\label{prop:torus}
Let $k>0$ and $l\ge 2$ be integers and suppose $0\le r\le \tfrac{1}{2}$.

(a) If $r\ne \tnfrac{2d+1}{2(2k+1)}$ for any integer $d$,
\[
\sigma_r(T_{2,2k+1})=
-2\left[ r(2k+1)+\tfrac{1}{2} \right].
\]

(b) 
For any integer $t$ with $1\le t\le \ttfrac{l}{2}$,
$\sigma_r
\left(T_{l,-l-1}\right)$ increases from $-2(t-1)^2+2l(t-1)$ to
$-2t^2+2(l+1)t-2$
over the interval $\tnfrac{t-1}{l}\le r \le \tdfrac{t}{l+1}$ and
decreases from
$-2t^2+2(l+1)t-2$ to $-2t^2+2lt$
over the interval $\tdfrac{t}{l+1}\le r\le \ttfrac{t}{l}$.
In particular, if
$\tnfrac{t-1}{l}\le r\le \ttfrac{t}{l}$, then
\[
-2(t-1)^2+2l(t-1) \le \sigma_r\left(T_{l,-l-1}\right)
\le -2t^2+2(l+1)t-2.
\]
\end{prop}

\begin{proof*}
Define
$f_{m,n}(r) = \frac{1}{2} \left(\text{jump in }
\sigma_r(T_{m,n}) \text{ at } r\right)$.  Then
\cite{lit79} showed that if $m,n>0$ and $0\le
r\le\tfrac{1}{2}$,
\[
f_{m,n}(r) =
\begin{cases}
(-1)^{\left[\ttfrac{a}{n}\right]+\left[\ttfrac{b}{m}\right]} &
\text{if } a,b,mnr \in\bbz,\ mr,nr \not\in\bbz 
\text{ with } mnr = am+bn, \\
0 & \text{otherwise}.
\end{cases}
\]
Note that $f_{m,n}(r)$ is nonzero only if $r=\ttfrac{s}{mn}$ for an
integer $s$ with $m\nmid s$ and $n\nmid s$.

To prove (a), let $m=2$ and $n=2k+1$.  Then
$f_{2,2k+1}(r)$ is nonzero only when $r=\ttfrac{s}{2(2k+1)}$ for odd $s$
with $1\le s\le 2k-1$. Note that $mnr=s = (-ks) (2) + (s)(2k+1)$.  So we
have
\[
f_{2,2k+1}\left(\frac{s}{2(2k+1)}\right) = 
(-1)^{[-ks/(2k+1)] + [s/2]}
\]
Write $s=2d-1$ for $1\le d\le k$.
Then
$[-ks/(2k+1)]=-d$ and $[s/2]=d-1$. Thus
\[
f_{2,2k+1}(r) = 
\begin{cases}
-1 & \text{if } r=\tnfrac{2d-1}{2(2k+1)} \text{ with } 1\le d\le k, \\
0 & \text{otherwise}.
\end{cases}
\]
Since $\tnfrac{2d-1}{2(2k+1)} \le r$ if and only if
$d \le r(2k+1)+\tfrac{1}{2}$, we see
\[
\sigma_r(T_{2,2k+1})=\sum_{d=1}^{\left[r(2k+1)+\ttfrac{1}{2}\right]}
2f_{2,2k+1}\left(\frac{2d-1}{2(2k+1)}\right) = 
-2\left[ r(2k+1)+\tfrac{1}{2}\right].
\]

To prove (b), 
note first that $f_{l,-l-1}=-f_{l,l+1}$ since $T_{l,-l-1}$ is the mirror
image of $T_{l,l+1}$.
Let $m=l$ and $n=l+1$.
For any integer $s$ with $0<s<\ttfrac{l(l+1)}{2}$, $l\nmid s$, and
$l+1\nmid s$, we see $s = (-s)(l)+(s)(l+1)$ and hence
\[
f_{l,-l-1}\left(\frac{s}{l(l+1)}\right) =
-f_{l,l+1}\left(\frac{s}{l(l+1)}\right) = -(-1)^{[-s/(l+1)]+[s/l]}.
\]
Let $t$ be an integer with $1\le t\le \ttfrac{l}{2}$.
Note $(l+1)(t-1)< lt < (l+1)t$.  We will consider two cases:
If $(l+1)(t-1)<s<lt$, then
$[-s/(l+1)]=-t$ and $[s/l]=t-1$ and hence
$f_{l,-l-1}\left(\ttfrac{s}{l(l+1)}\right) = 1$.
If $lt<s<(l+1)t$, then
$[-s/(l+1)]=-t$ and $[s/l]=t$ and hence
$f_{l,-l-1}\left(\ttfrac{s}{l(l+1)}\right) = -1$.

Since 
the sets $\{s\in\bbz\mid\tnfrac{t-1}{l} < \ttfrac{s}{l(l+1)}
<\tdfrac{t}{l+1}\}$ and 
$\{s\in\bbz\mid \tdfrac{t}{l+1} < \ttfrac{s}{l(l+1)} < \ttfrac{t}{l}\}$
have $l-t$ and $t-1$ elements, respectively,
$\sigma_r(T_{l,-l-1})$ changes by $2(l-t)-2(t-1)=2l-4t+2$ over the
interval $\tnfrac{t-1}{l}<r<\ttfrac{t}{l}$.  Incorporating these, we have
\begin{align*}
\sigma_r(T_{l,-l-1}) 
&=
\begin{cases}
\displaystyle
\sum_{t'=1}^{t-1} (2l-4t'+2)
& \text{if } r=\tnfrac{t-1}{l} \\
\displaystyle
\sum_{t'=1}^{t-1} (2l-4t'+2) + 2(l-t)
& \text{if } r=\ttfrac{t}{l+1}
\end{cases} \\
&=
\begin{cases}
-2(t-1)^2+2l(t-1) & \text{if } r=\tnfrac{t-1}{l} \\
-2t^2+2(l+1)t-2 & \text{if } r=\tdfrac{t}{l+1}.
\end{cases}
\end{align*}
Now (b) follows.
\end{proof*}

%%%%%%%%%%%%%%%%%% Section %%%%%%%%%%%%%%%%%%%%%%%%%

\section{Algebraically slice twisted doubles}
\label{sec:a-slice}
In this section we will estimate the Casson--Gordon invariants of
algebraically slice $D_k(K)$ and prove Theorem~\ref{thm:infinite},
Theorem~\ref{thm:independent}, and Corollary~\ref{cor:independent}.

As mentioned in Section~\ref{sec:double},
$D_k(K)$ is algebraically slice if
and only if $k=l(l+1)$ for an integer $l\ge 0$.
A Seifert matrix for $D_{l(l+1)}(K)$ corresponding to the Seifert surface
$F$ in Figure~\ref{fig:double} is
$\left(\begin{smallmatrix}
-1 & 1 \cr
0 & l(l+1)
\end{smallmatrix}\right)$.
The matrix $G=(A-A^t)^{-1}A$ associated with the isometric structure
$s$ is
$\left(\begin{smallmatrix}
0 & -l(l+1) \\ -1 & 1
\end{smallmatrix}\right)$ and has eigenvectors
\begin{align*}
v^+ &=\begin{pmatrix} l+1 \\ 1 \end{pmatrix}
\text{ corresponding to eigenvalue } -l \\
v^- &=\begin{pmatrix} -l \\ 1 \end{pmatrix}
\text{ corresponding to eigenvalue } l+1.
\end{align*}
With rational coefficients, $G$
is diagonalizable with respect
to the basis $\{v^+,v^-\}$.

For a positive integer $n$, let $nD_{l(l+1)}(K)$ be
the connected sum of $n$ copies of $D_{l(l+1)}(K)$.
We put a subscript $n$ on objects corresponding to  $nD_{l(l+1)}(K)$.
For example, $F_n$ denotes the Seifert surface of
$nD_{l(l+1)}(K)$ obtained by boundary connected summing $n$ copies
of the Seifert surface $F$ of $D_{l(l+1)}(K)$.

%%%%%%%%%%%%%% Subsection %%%%%%%%%%%%%%%%%%%%%

\subsection{Metabolizer of $nD_{l(l+1)}(K)$}
Let $Z_n$ be a metabolizer of the
associated isometric structure of $nD_{l(l+1)}(K)$.
(It should be remarked that unlike $F_n$, $s_n$, and $\theta_n$, $Z_n$
needs not be a direct sum of metabolizers of $D_{l(l+1)}(K)$.)
Since $G$ is diagonalizable with respect to the basis
$\{v^+,v^-\}$ with rational coefficients, $G_n$ associated
with $s_n$ is diagonalizable with respect to the basis
$\{v^+_j,v^-_j\}_{j=1,\ldots,n}$, where 
$
v^\pm_j=0\oplus \cdots \oplus v^\pm \oplus \cdots \oplus 0
$
is an eigenvector of $G_n$ whose $j$-th coordinate is the only nonzero
$v^\pm$ under the identification of $H_1(F_n)$ with the direct sum of $n$
copies of $H_1(F)$. Note that $\{v^+_1,\ldots,v^+_n\}$
is a basis for the eigenspace of $G_n$ corresponding to
the eigenvalue $-l$ and $\{v^-_1,\ldots,v^-_n\}$
is a basis for the eigenspace of $G_n$ corresponding to $l+1$.

The following lemma shows that $Z_n$ contains an eigenvector of $G_n$
for which the Casson--Gordon invariant can be easily estimated as will be
shown later.
\begin{lem}
\label{lem:metabolizer}
There are integers $e\ge\frac{n}{2}$, $a>0$, and
$a_{e+1},\ldots,a_n\in\bbz$ such that
$Z_n$ contains either
$a(\sum_{j=1}^e v^+_j) + \sum_{j=e+1}^n a_jv^+_j$
or $a(\sum_{j=1}^e v^-_j) + \sum_{j=e+1}^n a_jv^-_j$.
\end{lem}
\begin{proof*}
From a basic result of linear algebra,
$Z_n\otimes\bbq$ has a basis consisting of eigenvectors of $G_n$
since $Z_n$ is invariant under $G_n$ and $G_n$ is diagonalizable
over $\bbq$.
In particular, $Z_n\otimes\bbq = E^+\oplus E^-$, where $E^\pm$ are
the eigenspaces of $G_n$ restricted to $Z_n\otimes\bbq$.
Since the rank of $Z_n$ is $n$, one of $E^\pm$ has rank greater
than or equal to $\ttfrac{n}{2}$.
Suppose that $E^+$ has rank $e \ge \ttfrac{n}{2}$.
Using the Gauss-Jordan algorithm and rearranging basis elements, we may
assume that a basis for $E^+$ consists of vectors of the form
$v^+_j+u_j$, $1\le j\le e$, where $u_j$ are linear combinations of
$v^+_{e+1},\ldots,v^+_n$.
Adding these basis elements together we see that
$Z_n\otimes\bbq/\bbz$ contains a vector
$\sum_{j=1}^e v^+_j + \sum_{j=e+1}^n b_jv^+_j$ for some
$b_{e+1},\ldots,b_n\in\bbq$.
Multiplying the vector by a nonzero integer gives a desired vector.
The same argument works if rank$\,E^- \ge \ttfrac{n}{2}$.
\end{proof*}

%%%%%%%%%%%%%% Subsection %%%%%%%%%%%%%%%%%%%%%

\subsection{Estimation of the Casson--Gordon invariants}
We use Theorem~\ref{thm:naik}(b) to estimate $\tau$ of $D_{l(l+1)}(K)$ for
the character corresponding to $v^\pm$ in this subsection.

\begin{figure}
\vspace{1em}
\begin{tabular}{cc}
\setlength{\unitlength}{0.36pt}
\begin{picture}(464,180)
\put(-4,-2.5){\includegraphics[scale=0.36]{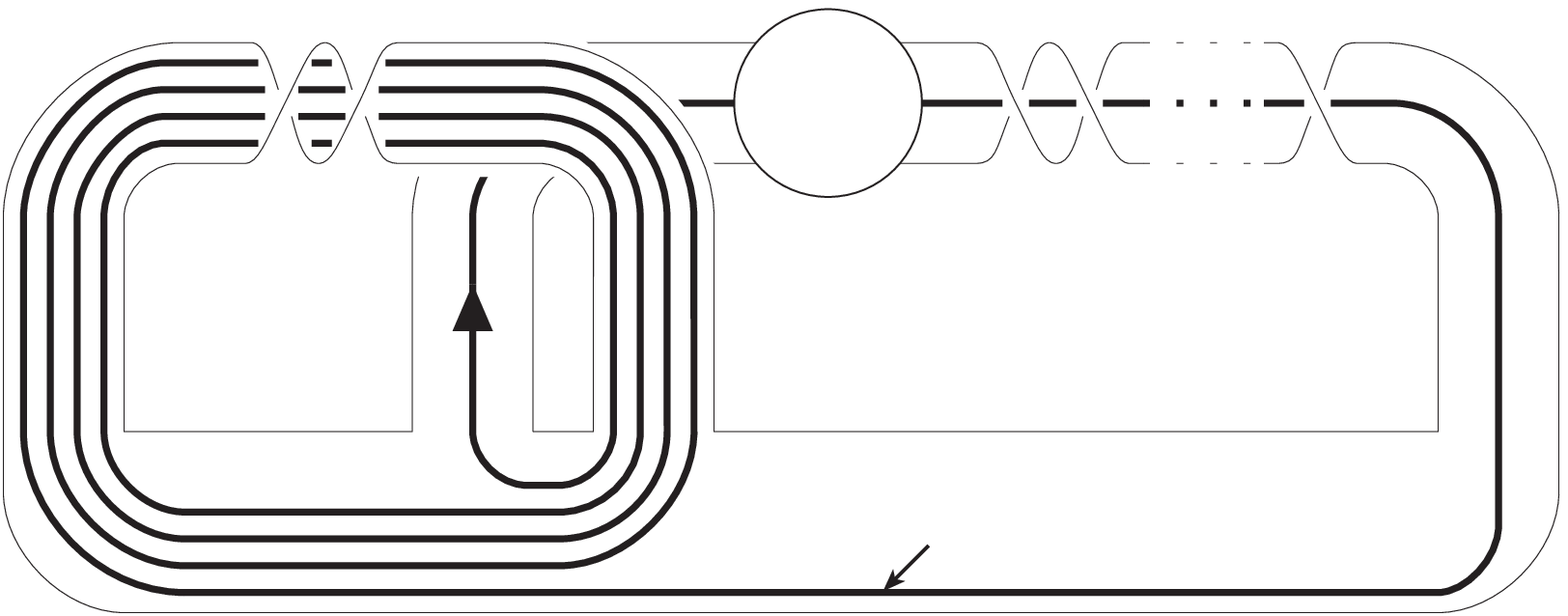}}
\put(342,122){\makebox(0,0)[t]{\tiny $k$ full twists}}
\put(246,152){\makebox(0,0){$K$}}
\put(276,20){\scriptsize $J_{x^+}$}
\end{picture}
&
\setlength{\unitlength}{0.36pt}
\begin{picture}(464,180)
\put(-4,-2.5){\includegraphics[scale=0.36]{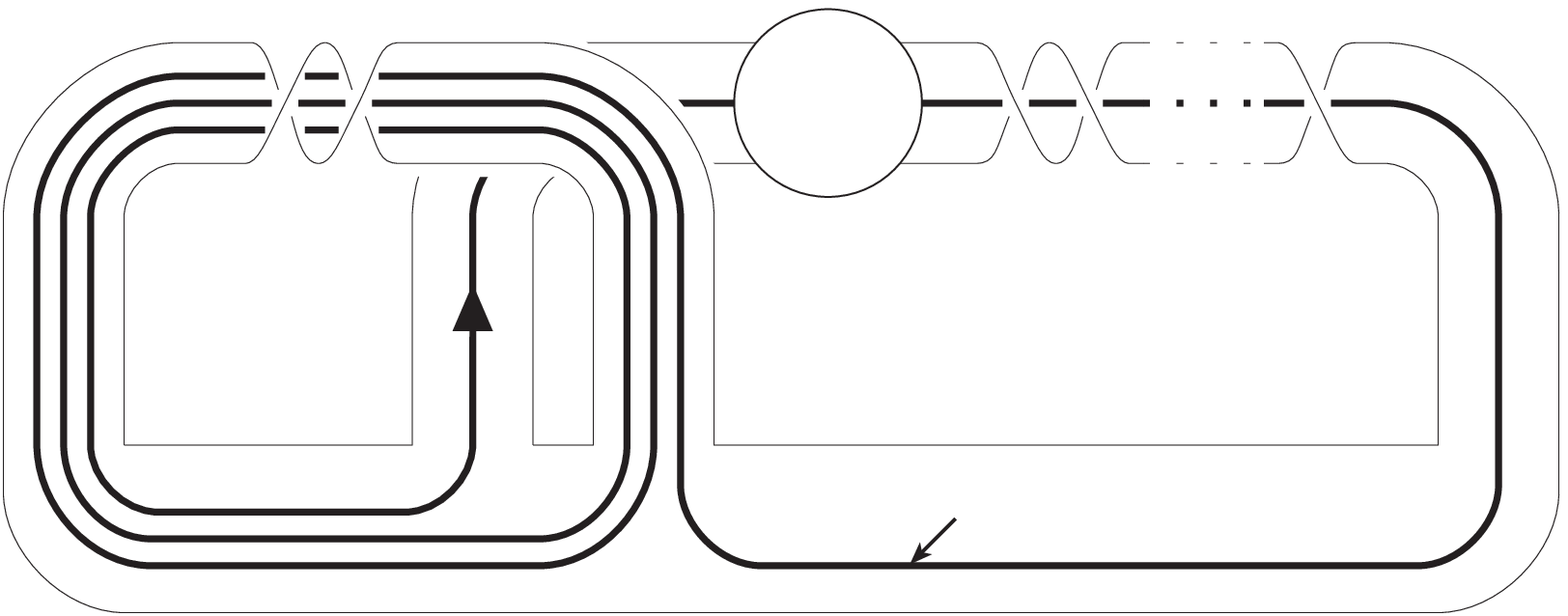}}
\put(342,122){\makebox(0,0)[t]{\tiny $k$ full twists}}
\put(246,152){\makebox(0,0){$K$}}
\put(284,28){\scriptsize $J_{x^-}$}
\end{picture}
\\[1ex]
$J_{x^+}$ for $x^+=(l+1,1)$ & $J_{x^-}$ for $x^-=(-l,1)$
\end{tabular}
\\[1em]
\setlength{\unitlength}{0.4pt}
\begin{picture}(306,144)
\put(-4,-2.5){\includegraphics[scale=0.4]{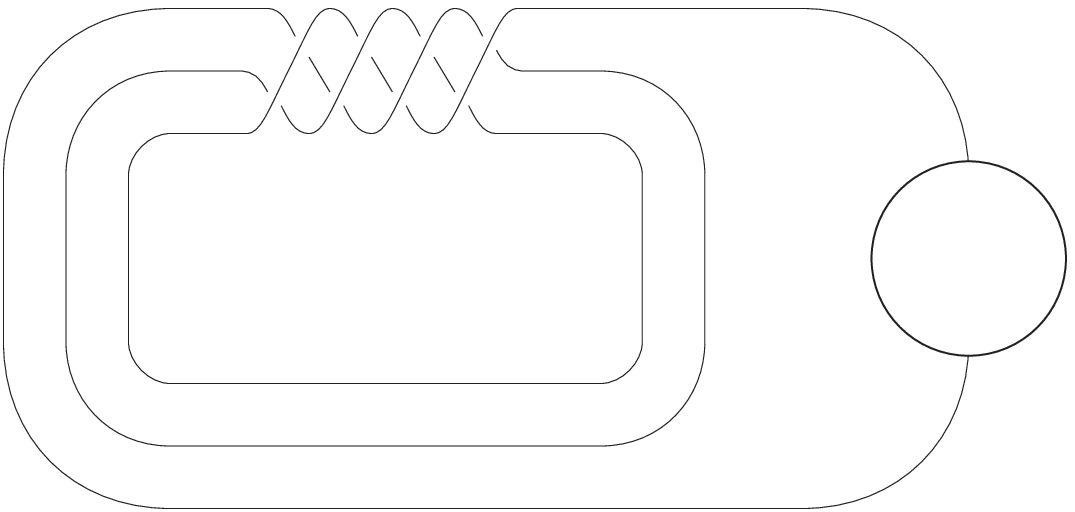}}
\put(278,72){\makebox(0,0){$K$}}
\end{picture}
\\[1ex]
$J_{x^\pm} = T_{l,-l-1}\# K$
\vspace{1em}
\caption{A knot $J_{x^\pm}$ for $l=3$}
\label{fig:alg-J_x}
\end{figure}

To change basis of $H_1(F)$ to either $\{x^+=v^+,y=(-1,0)\}$ or
$\{x^-=v^-,y=(-1,0)\}$, let
\[
P^+=\begin{pmatrix} l+1& -1 \\ 1 &0 \end{pmatrix}
\qquad \text{and} \qquad
P^-=\begin{pmatrix} -l & -1 \\ 1 &0 \end{pmatrix}.
\]
The Seifert matrices with respect to these bases $\{x^\pm,y\}$ are
\[
A^+=(P^+)^tAP^+ =
\begin{pmatrix} 0&l+1 \\ l&-1 \end{pmatrix}
\text{ and }
A^-=(P^-)^tAP^- =
\begin{pmatrix} 0&-l \\ -l-1&-1 \end{pmatrix},
\]
respectively.
Let $m^+=-l-1$ and $m^-=l$.
We can apply Theorem~\ref{thm:naik}(b) to
$A^\pm$ with $x=x^\pm$ and $m=m^\pm$, respectively, since $a=0$.  Observe
that we can choose
\[
J_{x^+} = J_{x^-} = T_{l,-l-1} \# K,
\]
as shown in Figure~\ref{fig:alg-J_x} when $l=3$.
Here, the property $T_{m,n}=T_{-m,-n}$ has been used.

%%%%%%%%%%%%%%%%%% Figure %%%%%%%%%%%%%%%%%%%%%%

For a prime power $q$, a prime $p$ dividing
$|(m+1)^q-m^q|=(l+1)^q-l^q$, and any integer $s$ relatively prime to $p$,
Theorem~\ref{thm:naik}(b) implies that
$x^\pm\otimes\ttfrac{s}{p}$ defines a character.  Abusing notation,
$x^\pm\otimes\ttfrac{s}{p}$ will also denote its corresponding character. 
We have
\[
\sigma_1\tau
\left(D_{l(l+1)}(K),x^\pm\otimes\frac{s}{p}\right) =
\sum_{i=0}^{q-1} \sigma_\frac{s_i}{p}
\left(T_{l,-l-1}\# K\right),
\]
where $s_i$ are integers such that
$0<s_i<p$ and $s_i\equiv (1+(m^\pm)^\ast)^i s$ \textup{mod} $p$.

Since a connected sum $K_1\# K_2$
can be considered as a satellite of $K_1$ with orbit $K_2$
and winding number $1$,  by Theorem~\ref{thm:satellite} we have

\begin{lem}
\label{lem:add}
For any prime power $q$, any prime $p$ dividing
$(l+1)^q-l^q$, and any integer $s$ with $p\nmid s$,
\[
\sigma_1\tau
\left(D_{l(l+1)}(K),x^\pm\otimes\frac{s}{p}\right) =
\sum_{i=0}^{q-1}\left(
\sigma_\frac{s_i}{p}
\left(T_{l,-l-1}\right)+
\sigma_\frac{s_i}{p}\left(K\right)
\right),
\]
where $s_i$ are integers such that
$0<s_i<p$ and $s_i\equiv (1+(m^\pm)^\ast)^i s$ \textup{mod} $p$.
\end{lem}

Now we need to estimate $\sigma_r(T_{l,-l-1})$.

\begin{lem}
\label{lem:estimate}
(a) 
\vspace{-1.3em}
\[
\sigma_r(T_{l,-l-1}) \ge
\begin{cases}
0 & \text{if } l\ge 1, \\
2 & \text{if } l\ge 2 \text{ and } \ttfrac{1}{l(l+1)}< r\le\tfrac{1}{2}.
\end{cases}
\]
(b)
For any constant $C_0$, there is an integer $l_0\ge 2$ such that,
for any $l\ge l_0$, any prime power~$q$, any prime $p$
dividing $(l+1)^q-l^q$, and any $m=m^\pm$, there is an
integer $s$ such that
\[
\sum_{i=0}^{q-1}\sigma_\frac{s_i}{p}(T_{l,-l-1}) > q C_0,
\]
where $s_i$ are integers such that
$0<s_i<p$ and $s_i\equiv (1+m^\ast)^is$ \textup{mod} $p$.
\end{lem}

\begin{proof*}
If $l=1$ then $T_{1,-2}$ is the unknot and $\sigma_r(T_{1,-2})=0$ for any
$r$.  Now suppose that $l\ge 2$.
From Proposition~\ref{prop:torus}(b)
$\sigma_r\left(T_{l,-l-1}\right)$ has local minima
$-2t^2+2lt$ at the integers $t$ with
$0\le t\le \ttfrac{l}{2}$.
Observe that the function
$-2t^2+2lt$ is increasing over $0\le t\le\ttfrac{l}{2}$
and has $0$ at $t=0$. 
Also, a close look at the proof of Proposition~\ref{prop:torus}(b) reveals
$f_{l,-l-1}(\ttfrac{1}{l(l+1)})=1$ 
and hence $\sigma_r(T_{l,-l-1})\ge 2$ if
$\ttfrac{1}{l(l+1)}< r\le \tfrac{1}{2}$.  Now (a) follows.

To show (b),
let $l_0$ be an integer such that $l_0\ge 2$ and
$\tfrac{3}{8}l_0^2 - l_0 -2  \ge 2 C_0$.
Let $l\ge l_0$, $q$ a prime power, and $p$ a prime dividing
$(l+1)^q-l^q$. 
It is easy to see that
$1+(m^\pm)^\ast \not\equiv 0$ \textup{mod} $p$.
For simplicity, let $a=1+(m^\pm)^\ast$
and $e$ the multiplicative order of $a$ \textup{mod} $p$.
Then $e$ divides $p-1$ and let $f=\tnfrac{p-1}{e}$.
For any integer $z$ not divisible by $p$, let
$E(z)=\{a^i z\in\bbz/p\}_{i=0,\ldots,e-1}$, that is, a coset of the
multiplicative group
$\bbz/p^\ast=\bbz/p-\{0\}$ modulo $\langle a\rangle$. Then there are
integers $z_1,\ldots,z_f$ such that
$\cup_jE(z_j) = \bbz/p^\ast$ and $E(z_j)$ are all disjoint.  

Let $P=\{x\in\bbz/p^\ast\mid \ttfrac{p}{4}\le x\le \ttfrac{3p}{4}\}$.
We will show that there is $s$ for which at least half of
$s_0,\ldots,s_{q-1}$ belong to $P$.  Note $|P|\ge
\tnfrac{p-1}{2}=\ttfrac{|\bbz/p^\ast|}{2}$. Since
$\{E(z_j)\}_{j=1,\ldots,f}$ is a partition of $\bbz/p^\ast$, there is
$j$ such that $|E(z_j)\cap P|\ge\ttfrac{|E(z_j)|}{2} = \ttfrac{e}{2}$.
For simplicity, assume $j=1$.

Let $G=\{(c,d)\in\bbz\times\bbz \mid 0\le c\le q-1 \text{ and } 0\le d\le
e-1\}$.  Define a function $\phi\!:G\to E(z_1)$ by $\phi(c,d)=a^ca^dz_1$
\textup{mod} $p$.  Since, for each fixed $c$,
$\{\phi(c,d)\mid 0\le d\le e-1\}=E(z_1)$, the function $\phi$ is a $q$ to
$1$ map.  Since $|E(z_1)\cap P|\ge\ttfrac{|E(z_1)|}{2}$,
$|\phi^{-1}(E(z_1)\cap P)| \ge 
q\left(\ttfrac{|E(z_1)|}{2}\right)=\ttfrac{|G|}{2}$.

For each fixed $d$ with $0\le d\le e-1$, let
$G_d=\{(c,d)\in G\mid 0\le c\le q-1\}$.
Then $\{G_d\}_{d=0,\ldots,e-1}$ is a partition of $G$ and hence there is
an integer $d_0$ such that
$|G_{d_0}\cap \phi^{-1}(E(z_1)\cap P)| \ge
\ttfrac{|G_{d_0}|}{2}=\ttfrac{q}{2}$.
Let $s=a^{d_0}z_1$ and, for $i=0,\ldots,q-1$, let $s_i$ be integers
such that $0<s_i<p$ and $s_i\equiv a^i s$ \textup{mod} $p$.  Then
$\phi(G_{d_0})=\{s_0,s_1,\ldots,s_{q-1}\}\subset\bbz/p$ and at least half
of $s_i$'s belong to $P$.

Note that if $s_i\in P$, then $\tfrac{1}{4}\le \ttfrac{s_i}{p}\le
\tfrac{3}{4}$.  Let $t=\left[\ttfrac{l}{4}+1\right]$.
Then $\ttfrac{l}{4}<t\le\ttfrac{l}{4}+1$ and 
$\tnfrac{t-1}{l}\le \tfrac{1}{4}\le \ttfrac{s_i}{p}$.
By Proposition~\ref{prop:torus}(b), if $s_i\in P$,
\begin{align*}
\sigma_\frac{s_i}{p}(T_{l,-l-1}) &
\ge -2(t-1)^2+2l(t-1) \\
&
> -2\left(l/4-1\right)^2 + 2l\left(l/4-1\right) 
&& \text{since } l/4 < t\\
& = \tfrac{3}{8}l^2 - l -2 \\
& \ge 2C_0
&& \text{by the definition of } l_0.
\end{align*}
Summing these, we have
\begin{align*}
\sum_{i=0}^{q-1} \sigma_\frac{s_i}{p}(T_{l,-l-1}) &
= \left(\sum_{s_i\in P} + \sum_{s_i\not\in P}\right)
\sigma_\frac{s_i}{p}(T_{l,-l-1}) \\
&
\ge \sum_{s_i\in P}  \sigma_\frac{s_i}{p}(T_{l,-l-1}) 
&& \quad\text{by (a)}\\
&
> \frac{q}{2}\cdot 2C_0=qC_0.
\end{align*}
This completes the proof.
\end{proof*}

%%%%%%%%%%%% Subsection %%%%%%%%%%%%%%%%%%

\subsection{Homology of prime power
fold cyclic branched covers}
We need some algebraic background.
The \emph{resultant} of two non-constant integral polynomials $f(t)$ and
$g(t)$ is defined as follows:  We may factor completely the polynomials
$f$ and $g$ in some extension ring of $\bbz$ as:
$f(t)=a\prod^n_{i=1}(t-\alpha_i)$ and $g(t)=b\prod^m_{j=1}(t-\beta_j)$.
Then the \emph{resultant} of $f$ and $g$, denoted $R(f,g)$, is
$a^mb^n\prod^n_{i=1}\prod^m_{j=1}(\alpha_i-\beta_j)$.
It is easy to see that $R(f,g)=0$ if and only if $f$ and $g$
have a common root in a field over $\bbz$.  We remark that this is also
valid when working modulo a prime $p$.

It is known by Fox (see~\cite{web79} for a proof)
that the order of the homology of the $n$-fold cyclic cover of
$S^3$ branched over a knot $J$ is the absolute value of the resultant,
$|R(t^n-1,\Delta_J(t))|$, of $t^n-1$ and the Alexander polynomial
$\Delta_J(t)$ of $J$.

\begin{prop}
\label{prop:larsen}
For any knot $K$ and any integer $k>0$, there are infinitely many distinct
primes each of which divides $|H_1(M^q)|$ for some prime $q$, where $M^q$
is the $q$-fold cyclic cover of $S^3$ branched over $D_k(K)$.
\end{prop}

\begin{proof*}
The Alexander polynomial of $D_k(K)$ is $\Delta_k=-kt^2+(2k+1)t-k$.
Let $R_n$ denote $|R(t^n-1,\Delta_k)|$.

First, we will show that $R_{q_1}$ and $R_{q_2}$ are relatively
prime for distinct primes $q_1$ and $q_2$.  
Suppose to the contrary there is a prime $p$ dividing both $R_{q_1}$ and
$R_{q_2}$. Then, working modulo $p$, for $j=1,2$, $\Delta_k$ and
$t^{q_j}-1$ have a common root $r_j$ in an extension field of the
finite field $\mathbf{F}_p$ of $p$ elements.

We claim that the three polynomial $\Delta_k$,
$t^{q_1}-1$, and $t^{q_2}-1$ have a common root, $r$, \textup{mod} $p$. 
If
$r_1=r_2$, this is obvious.  Suppose $r_1\ne r_2$. Then $\Delta_k$ has
two distinct roots $r_1$ and $r_2$ \textup{mod} $p$ and hence
$\Delta_k$ must be quadratic over $\mathbf{F}_p$.  In particular, $k\ne
0$ in $\mathbf{F}_p$.  Thus, $r_1r_2=\ndfrac{-k}{-k}=1$ or
$r_2=\ttfrac{1}{r_1}$. So, $r_2^{q_1}=\ttfrac{1}{r_1^{q_1}}=1$
\textup{mod} $p$ and
$r_2$ is a common root, $r$, \textup{mod} $p$ of the three polynomials.

Since $q_1$ and $q_2$ are distinct primes, there are 
$a$ and $b$ such that $aq_1+bq_2=1$.  So
$r=r^{aq_1+bq_2}=(r^{q_1})^a(r^{q_2})^b\equiv 1$ \textup{mod} $p$.
This implies that $1$ is a root of $\Delta_k$ mod~$p$.
However, this is a contradiction since
$\Delta_k(1)=1\not\equiv 0$ \textup{mod} $p$.
Thus there are no primes $p$ dividing both $R_{q_1}$ and $R_{q_2}$,
implying they are relatively prime.

It now suffices to show that $R_q>1$ for any large primes $q$.  We
will show that $R_q\to\infty$ as prime
$q\to\infty$.  If $G$ is a matrix associated with the isometric
structure of $D_k(K)$ and $N^q$ is the kernel of
$(G^q-(G-I)^q)\otimes \id_{\bbq/\bbz}$ as before,
then $R_q=|\det\left(G^q-(G-I)^q\right)|$ since
$R_q=|H_1(M^q)|=|N^q|$. A Seifert matrix for $D_k(K)$
corresponding to the Seifert surface in Figure~\ref{fig:double} is
$A =\left(\begin{smallmatrix}
-1 & 1 \cr 0 & k
\end{smallmatrix}\right).$
So $G=(A-A^t)^{-1}A=\left(\begin{smallmatrix}
0 & -k \cr -1 & 1
\end{smallmatrix}\right).$
Let $u=\tnfrac{1+\sqrt{4k+1}}{2}$, $w=\tnfrac{1-\sqrt{4k+1}}{2}$,
and
$P =\left(\begin{smallmatrix}
w & u \cr 1 & 1
\end{smallmatrix}\right).$
Then $P^{-1}GP =\left(\begin{smallmatrix}
u & 0 \\ 0 & w \end{smallmatrix}\right)$
and
$P^{-1}(G-I)P =\left(\begin{smallmatrix}
-w & 0 \\ 0 & -u \end{smallmatrix}\right)$.
We now see that, for any odd integer $q$, 
\[
P^{-1}(G^q-(G-I)^q)P = \left(u^q+w^q\right)I
\]
and hence
$R_q = (u^q+w^q)^2$.
Since
$|u|>1$, $\left|\ttfrac{w}{u}\right| < 1$, and
\[
\sqrt{R_q} =
|u^q+w^q|\ge |u|^q-|w|^q=|u|^q\left(1-\left|\frac{w}{u}\right|^q\right),
\]
$R_q\to\infty$ as prime $q\to\infty$.
In particular, $R_q> 1$ for any large primes $q$.  This completes the
proof.
\end{proof*}

%%%%%%%%%%%%% Subsection %%%%%%%%%%%%%%%%

\subsection{Proofs of Theorem~\ref{thm:infinite},
Theorem~\ref{thm:independent}, and Corollary~\ref{cor:independent}}

\begin{proof*}[Proof of Theorem~\ref{thm:infinite}.]
Let $\scrm=\min_{0<r<1} \sigma_r(K)$.
For $C_0=2|\scrm|$, there is $l_0$ satisfying the conclusion of
Lemma~\ref{lem:estimate}(b).  Let $l\ge l_0$.

Suppose to the contrary
that $n>0$ and $nD_{l(l+1)}(K)$ is zero in $\Gamma^+$.
Then $nD_{l(l+1)}(K)$ has vanishing Gilmer
slice obstruction.  
Use the same notation $F$, $M^q$, $N$, etc.~as before for
$D_{l(l+1)}(K)$ and put a subscript $n$ on the objects corresponding to
$nD_{l(l+1)}(K)$. Let $Z_n$ be a metabolizer satisfying the conclusion of
Theorem~\ref{thm:gilmer}(a) for the surface $F_n$.
By Lemma~\ref{lem:metabolizer},
$Z_n$ contains an integral vector $v$ that is either
$a(\sum_{j=1}^e v^+_j) + \sum_{j=e+1}^n a_jv^+_j$
or $a(\sum_{j=1}^e v^-_j) + \sum_{j=e+1}^n a_jv^-_j$
for some $a> 0$ and $e\ge\ttfrac{n}{2}$.
By Proposition~\ref{prop:larsen} we can find a prime
$p$ and an odd prime $q$ such that $p$ divides $|H_1(M^q)|$,
$p$ does not divide $a$, and $\tau(nD_{l(l+1)}(K), (N^q_n)_p\cap
(Z_n\otimes\bbq/\bbz))$ vanishes, where $(N^q_n)_p$ denotes the
$p$-primary component subgroup of $N^q_n$.

Recall that $N^q \cong \ker ((G^q-(G-I)^q)\otimes \id_{\bbq/\bbz})$.
From the proof of Proposition~\ref{prop:larsen},
$G^q-(G-I)^q$ is the identity matrix
multiplied by an integer $h=|u^q+w^q|=(l+1)^q-l^q$ when $q$ is odd
and $k=l(l+1)$.
Thus $N^q=\{z\otimes \ttfrac{1}{h}\mid z\in H_1(F)\}$.
Since $N^q_n$ is the direct sum of $n$ copies of $N^q$, we have
$N^q_n=\{z_n\otimes \ttfrac{1}{h} \mid z_n\in H_1(F_n)\}$.
Since the prime $p$ divides $|H_1(M^q)|=h^2$, $p$ divides $h$ and hence
$z_n\otimes \ttfrac{1}{p}\in (N^q_n)_p$ for any $z_n\in H_1(F_n)$.
In particular, for the vector $v\in Z_n$ chosen above,
$v\otimes\ttfrac{1}{p}\in (N^q_n)_p$.
Moreover, since $v\in Z_n$,
$v\otimes\ttfrac{1}{p}\in (N^q_n)_p\cap (Z_n\otimes\bbq/\bbz)$ and hence
$\sigma_1\tau\left(nD_{l(l+1)}(K),v\otimes\ttfrac{1}{p}\right)=0$.

On the other hand,
let $s$ be the constant from Lemma~\ref{lem:estimate}(b) determined by
$l$, $p$, and $q$ chosen above together with $m=m^+$ or $m^-$ depending on
whether
$v=a(\sum_{j=1}^e v^+_j) + \sum_{j=e+1}^n a_jv^+_j$
or $a(\sum_{j=1}^e v^-_j) + \sum_{j=e+1}^n a_jv^-_j$. Since
$p$ was chosen not to divide $a$, by multiplying an integer to $v$, we
may further assume that $a\equiv s$ \textup{mod} $p$.
For $i=0,\ldots,q-1$, $j=e+1,\ldots,n$, let $s_i$, $s_{ij}$ be integers
such that $0<s_i, s_{ij}< p$, $s_i\equiv (1+m^\ast)^ia$ \textup{mod} $p$,
and $s_{ij}\equiv (1+m^\ast)^i a_j$ \textup{mod} $p$.
Since $v\otimes\ttfrac{1}{p}=\sum_{j=1}^e v^\pm_j\otimes\ttfrac{a}{p}
+ \sum_{j=e+1}^n v^\pm_j\otimes\ttfrac{a_j}{p}$,
we have
\[
\sigma_1\tau\left(nD_{l(l+1)}(K),v\otimes\frac{1}{p}\right)
= \sigma_1\tau\left(nD_{l(l+1)}(K),
\sum_{j=1}^e v^\pm_j\otimes\frac{a}{p}
+ \sum_{j=e+1}^n v^\pm_j\otimes\frac{a_j}{p}\right).
\]
By the additivity of $\sigma_1\tau$, it is equal to
\[
\sum_{j=1}^e \sigma_1\tau
\left(D_{l(l+1)}(K),x^\pm\otimes\frac{a}{p}\right) 
+
\sum_{j=e+1}^n \sigma_1\tau
\left(D_{l(l+1)}(K),x^\pm\otimes\frac{a_j}{p}\right),
\]
which is, by
Lemma~\ref{lem:add}, equal to
\[
\sum_{j=1}^e \left( \sum_{i=0}^{q-1}\sigma_\frac{s_i}{p}(T_{l,-l-1})
+ \sum_{i=0}^{q-1}\sigma_\frac{s_i}{p}(K)\right)
+ \sum_{j=e+1}^n \left(
\sum_{i=0}^{q-1}\sigma_\frac{s_{ij}}{p}(T_{l,-l-1})
+ \sum_{i=0}^{q-1}\sigma_\frac{s_{ij}}{p}(K) \right).
\]
By Lemma~\ref{lem:estimate}, we now see that
\[
\sigma_1\tau\left(nD_{l(l+1)}(K),v\otimes\frac{1}{p}\right)
> e q C_0 + nq \scrm 
\ge q( 2e |\scrm| + n \scrm)
\ge 0.
\]
Thus
$\sigma_1\tau\left(nD_{l(l+1)}(K),v\otimes\ttfrac{1}{p}\right)\ne 0$,
which is a contradiction. This proves~(a).

Next, suppose that $l\ge 2$ and $\sigma_r(K)\ge 0$ for all $r$.
Under the same
contradiction hypothesis and notation as above except:
let $s=\left[\ttfrac{p}{2}\right]$, instead of choosing it from
Lemma~\ref{lem:estimate}(b).  Observe that
$\ttfrac{1}{l(l+1)}<\tfrac{1}{3} \le \ttfrac{s}{p}\le \tfrac{1}{2}$ and
hence
$\sigma_{\ttfrac{s}{p}}(T_{l,-l-1})\ge 2$ by Lemma~\ref{lem:estimate}(a).
We may assume $a\equiv s$ \textup{mod} $p$ as before.  Then we have
\[
\sigma_1\tau\left(nD_{l(l+1)}(K),v\otimes\frac{1}{p}\right)
=
\sigma_\frac{s_0}{p}(T_{l,-l-1}) + \text{ other terms }
\ge
\sigma_\frac{s}{p}(T_{l,-l-1})
\ge 2.
\]
This is a contradiction, 
completing the proof of the first part of (b).

In addition, suppose that $\sigma_{r_0}(K)>0$ for some $r_0$.
We only need to check the case $l=1$.
By the definition of the averaged signature $\sigma_r$,
there are $r_1$ and $r_2$ such that $0<r_1<r_2\le\tfrac{1}{2}$ and
$\sigma_r(K)>0$ for any $r$ with $r_1 <r<r_2$.
We can take $p$ arbitrary large as above so that $r_1
<\ttfrac{s}{p}<r_2$ for some integer $s$.  Apply a similar argument as
above.  The difference from the previous case is that the role of
$\sigma_{\ttfrac{s}{p}}(T_{l,-l-1})$ is switched with that of
$\sigma_{\ttfrac{s}{p}}(K)$.
\end{proof*}

\begin{proof*}[Proof of Theorem~\ref{thm:independent}.]
The Alexander polynomial of $D_k(K)$ is $\Delta_k = -kt^2+(2k+1)t-k$.  We
will see that all $\Delta_k$ are coprime in $\bbq[t^{\pm1}]$.
Let $k\ne l$ and let $g$ be the greatest common divisor
of $\Delta_k$ and $\Delta_l$ in $\bbq[t^{\pm1}]$.
Then $g$ divides $l\Delta_k-k\Delta_l=(l-k)t$ that is a unit in
$\bbq[t^{\pm1}]$.  So,
$\Delta_k$ and $\Delta_l$ are relatively prime in
$\bbq[t^{\pm1}]$ for any distinct pair of integers $k$ and $l$.

Let $\scrl$ be the set of integers $l\ge 0$ for which $D_{l(l+1)}(K)$ has
infinite order in $\Gamma^+$.
Theorem~\ref{thm:infinite} implies
$\scrl$ contains all but finitely many nonnegative
integers.  Note that under the hypotheses of (b) $\scrl$ contains
all but $0,1$ or only $0$.
Suppose that there are distinct $l_i\in\mathcal{L}$ such that
$e_1D_{l_1(l_1+1)}(K)\# e_2D_{l_2(l_2+1)}(K)\#\cdots
\#e_nD_{l_n(l_n+1)}(K)$ is zero in $\Gamma^+$ for some integers
$e_1,\ldots,e_n$.
Since $\Delta_{e_iD_{l_i(l_i+1)}(K)}(t)=
\left( \Delta_{l_i(l_i+1)}\right)^{e_i}$,
all $\Delta_{e_iD_{l_i(l_i+1)}(K)}(t)$ are pairwise coprime.
Also, note that all $e_iD_{l_i(l_i+1)}(K)$ have nonsingular Seifert forms.
Applying
Theorem~\ref{thm:alexander} inductively,
each $e_iD_{k_i}(K)$ must be zero in $\Gamma^+$
and hence each $e_i$ must be $0$.
This completes the proof.
\end{proof*}

\begin{proof*}[Proof of Corollary~\ref{cor:independent}.]
Part (a) is an immediate consequence of
Theorem~\ref{thm:independent}(b) since $\sigma_r(K)=0$ for any $r$
if $K$ is the unknot.  For part~(b), the torus knots $T_{l,-l-1}$
for $l>1$ are such examples with $\sigma_r \ge 0$ for all $r$ and
$\sigma_r>0$ for some $r$.
\end{proof*}

%%%%%%%%%%%%%%%%%% Section %%%%%%%%%%%%%%%%%%%%%%%%%

\section{Non-ribbonness of linear combinations of twisted doubles}
\label{sec:ribbon}
In this section we will estimate the Casson--Gordon invariants
of all twisted doubles $D_k(K)$ for double branched covers and prove
Theorem~\ref{thm:rinfinite} and Theorem~\ref{thm:ribbon}.

%%%%%%%%%%%%%%%%%% Subsection %%%%%%%%%%%%%%%%%%%%%%%%%%

\subsection{Estimation of the Casson--Gordon invariants}
\label{sec:rest}
We estimate the Casson--Gordon invariants of $D_k(K)$
for $k>0$ in this subsection.  Recall from Theorem~\ref{thm:levine}
that if $k\ge 0$, then $D_k(K)$ has finite order in the algebraic
concordance group~$\scrg$.

A Seifert matrix for $D_k(K)$ corresponding to the Seifert surface in
Figure~\ref{fig:double} is
$\left(\begin{smallmatrix}
-1 & 1 \cr
0 & k
\end{smallmatrix}\right).$
By changing basis to $\{x=(1,2),y=(0,1)\}$, the Seifert matrix changes
to the matrix
\[
A =
\left(\begin{matrix} 1 & 0 \cr 2 & 1 \end{matrix}\right)^t
\left(\begin{matrix} -1 & 1 \cr 0 & k \end{matrix}\right)
\left(\begin{matrix} 1 & 0 \cr 2 & 1 \end{matrix}\right)
=
\left(\begin{matrix} 4k+1 & 2k+1 \cr 2k & k \end{matrix}\right).
\]
In this case, $a=4k+1$, $m=-(2k+1)$, and $b=k$ following the notation of
Theorem~\ref{thm:naik}.

\begin{figure}
\vspace{1em}
\begin{tabular}{cc}
\setlength{\unitlength}{0.4pt}
\begin{picture}(464,180)
\put(-4,-2.5){\includegraphics[scale=0.4]{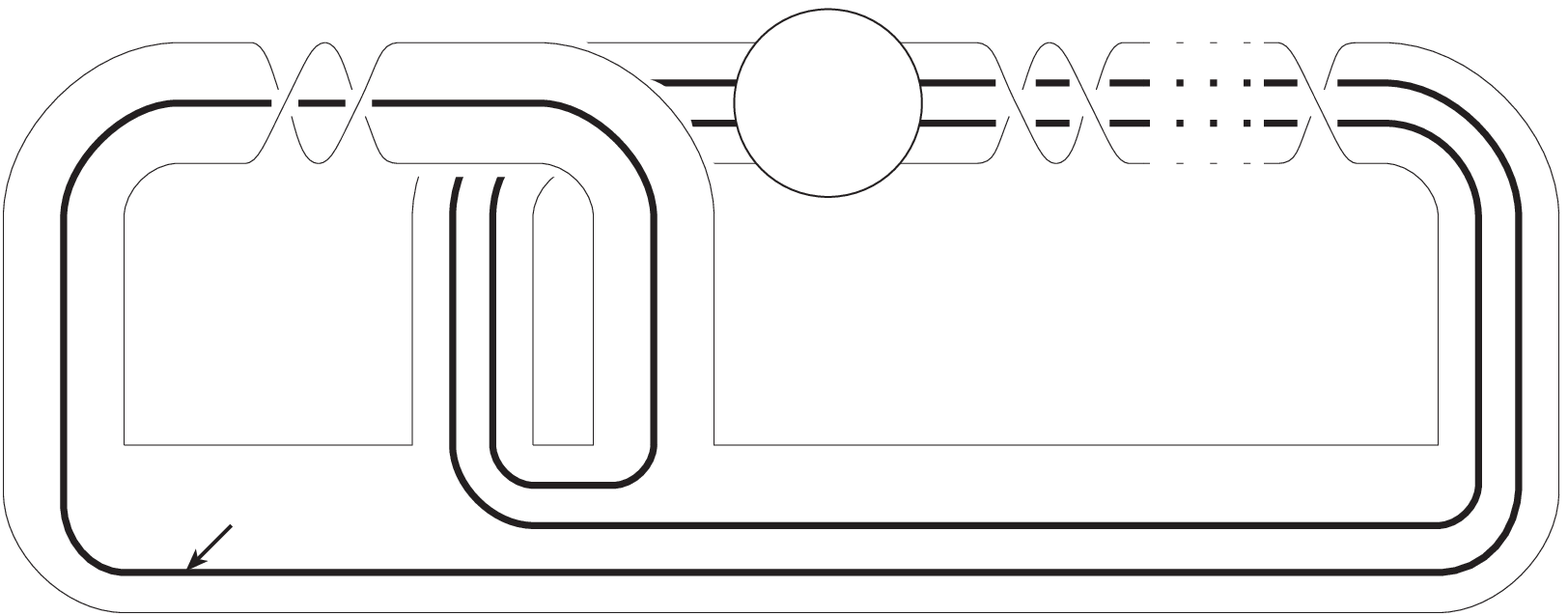}}
\put(342,122){\makebox(0,0)[t]{\tiny $k$ full twists}}
\put(246,152){\makebox(0,0){$K$}}
\put(68,26){\scriptsize $J_x$}
\end{picture}
&
\setlength{\unitlength}{0.4pt}
\begin{picture}(306,180)
\put(-4,-2.5){\includegraphics[scale=0.4]{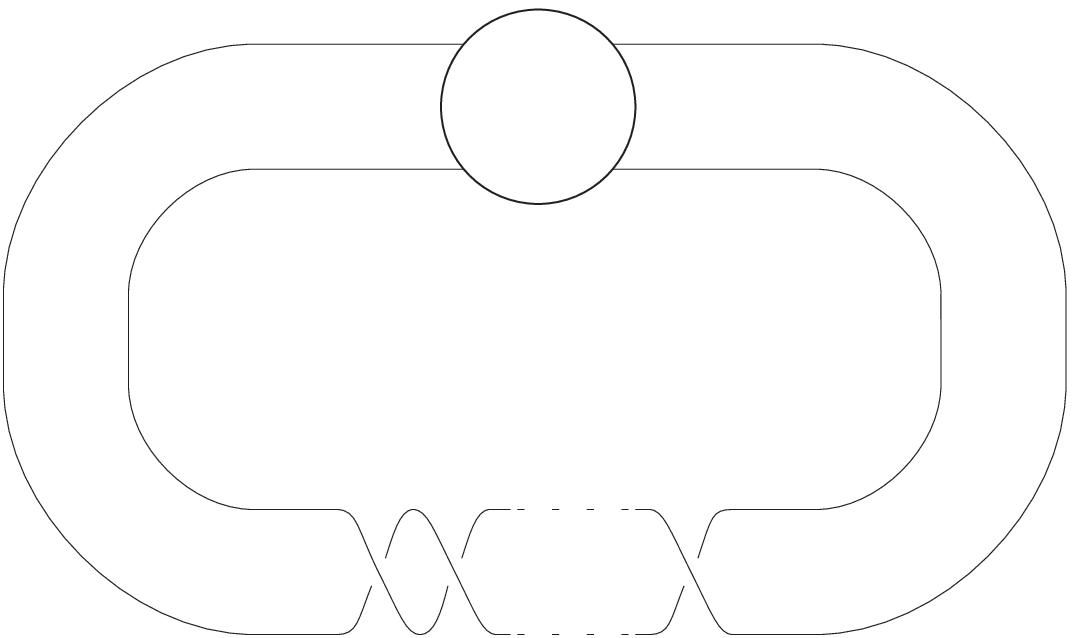}}
\put(154,45){\makebox(0,0)[b]{\tiny $k$ full twists}}
\put(154,152){\makebox(0,0){$K$}}
\end{picture}
\\[1ex]
$J_x$ in $D_k(K)$ & $J_x = K(T_{2,2k+1})$
\end{tabular}
\vspace{1em}
\caption{A knot $J_x$ that represents $x=(1,2)$}
\label{fig:J_x}
\end{figure}

We consider only the case $q=2$.
The map $\varepsilon^2$ is represented by the matrix
\[ 
G^2-(G-I)^2 =
\left[\left(A-A^t\right)^{-1}A\right]^2 -
\left[\left(A-A^t\right)^{-1}A^t\right]^2 =
\left(\begin{matrix} -(4k+1) & -2k \cr 2(4k+1) & 4k+1 \end{matrix}\right).
\]
Let $p$ be a prime dividing $4k+1$ and let $s$ be an integer with $0<s<p$.
Note that $p$ is odd and $\varepsilon^2\otimes \id_{\bbq/\bbz}
(x\otimes \ttfrac{s}{p})=0$
in $H_1(F)\otimes \bbq/\bbz$, i.e.~
$x\otimes \ttfrac{s}{p}$ is in the kernel of $\varepsilon^2\otimes
\id_{\bbq/\bbz}$. Note that $J_x$, a simple closed curve on $F$
representing $x=(1,2)$, can be represented by $K(T_{2,2k+1})$, a
satellite knot of $K$ with orbit
$T_{2,2k+1}$ as shown in Figure~\ref{fig:J_x}.
Note $\sigma_{\ttfrac{s}{p}}(K(T_{2,2k+1}))=\sigma_{\ttfrac{2s}{p}}(K)+
\sigma_{\ttfrac{s}{p}}(T_{2,2k+1})$ by Theorem~\ref{thm:satellite}.

%%%%%%%%%%%%%%%%% Figure %%%%%%%%%%%%%%%%%%%%%%

By Theorem~\ref{thm:naik}, $x\otimes\ttfrac{s}{p}$ defines a character
$\chi_{\ttfrac{s}{p}}$ and
\[
\sigma_1\tau(D_k(K),\chi_\frac{s}{p}) =
2\sigma_{\frac{2s}{p}}(K) +
2\sigma_\frac{s}{p}(T_{2,2k+1})+\frac{4(p-s)s(4k+1)}{p^2} -
\sigma_\frac{1}{2}(D_k(K)).
\]
Since $D_k(K)$ has finite algebraic order and since
$\sigma_{\frac{1}{2}}$ is additive under connected sums,
$\sigma_\frac{1}{2}(D_k(K)) = 0$.
Thus, we have
\[
\sigma_1\tau(D_k(K),\chi_\frac{s}{p})=
2\sigma_{\frac{2s}{p}}(K) +
2\sigma_\frac{s}{p}(T_{2,2k+1}) +
4\left(\frac{s}{p}\right)\left(1-\frac{s}{p}\right)(4k+1).
\]
We will show

\begin{lem}
\label{lem:minmax}
Let $\scrm = 2\min_{0<r<1} \sigma_r(K)$.

(a) For any $k\ge 3$,
\[
\min_\chi \sigma_1\tau(D_k(K),\chi) \ge \scrm -\frac{4}{4k+1},
\]
where $\chi$ runs over all prime power characters.

(b) Let $s$ be an integer such that $p=4s\pm 1$.
Then, for any constant $C_0$, there is $k_0\ge 1$ such that,
for any $k\ge k_0$,
\[
\sigma_1\tau(D_k(K),\chi_\frac{s}{p}) > C_0.
\]

(c) Suppose that $\scrm \ge 0$. If $k\ge 3$ then
$\sigma_1\tau(D_k(K),\chi_r)$ may have only one
non-positive value $-\tdfrac{4}{4k+1}$ at $r=\tdfrac{1}{4k+1}$, and
if we let $c=\tnfrac{4k+1}{p}$
\[
\sigma_1\tau(D_k(K),\chi_\frac{s}{p}) >
\begin{cases}
cs-2 &  \text{if } p=4s+1, \\[1ex]
cs-c/2-2 & \text{if } p=4s-1.
\end{cases}
\]
\end{lem}

\begin{proof*}
Let $r=\ttfrac{s}{p}$, where $0<s<p$.
Since $\sigma_1\tau(D_k(K),\chi_r)=\sigma_1\tau(D_k(K),\chi_{1-r})$, it
suffices to compute $\sigma_1\tau$ for $\chi_\ttfrac{s}{p}$ when
$1\le s\le \tnfrac{p-1}{2}$.
Then $\tdfrac{1}{4k+1}\le r \le\tdfrac{2k}{4k+1}$.
For $\tdfrac{1}{4k+1}\le r\le \tdfrac{2k}{4k+1}$, let
\[
f(r)
= \tfrac{1}{2}\left(\sigma_1\tau(D_k(K),\chi_r)-2\sigma_{2r}(K)\right)
= \sigma_r(T_{2,2k+1}) +2r(1-r)(4k+1).
\]
From Proposition~\ref{prop:torus}(a), we have
$-2(2k+1)r-1\le \sigma_r(T_{2,2k+1})$.
Let
\[
g(r) =-2(2k+1)r-1+2r(1-r)(4k+1)=-2(4k+1)r^2+4kr-1.
\]
Then $g(r)\le f(r)$.
Observe that $g$ is
a quadratic polynomial in $r$ with maximum at $r=\tdfrac{k}{4k+1}$ and
that
\[
g\left(\frac{s}{4k+1}\right) =
\begin{cases}
-1 & \text{if } s=2k, \\
-\tdfrac{3}{4k+1} & \text{if } s=1 \text{ and } 2k-1, \\
\ndfrac{4k-9}{4k+1} & \text{if } s= 2 \text{ and } 2k-2.
\end{cases}
\]

To prove (a) and the first part of (c), assume that $k\ge 3$.
Then $g\left(\tdfrac{2}{4k+1}\right)=g\left(\ndfrac{2k-2}{4k+1}\right)=
\ndfrac{4k-9}{4k+1} > 0$
and hence $f(r) \ge g(r) >0$ if $\tdfrac{2}{4k+1} \le r \le
\ndfrac{2k-2}{4k+1}$. Now, we will compute $f(r)$
when $r=\tdfrac{1}{4k+1}$, $\ndfrac{2k-1}{4k+1},$ and $\tdfrac{2k}{4k+1}$.
By Proposition~\ref{prop:torus}(a),
\[
\sigma_r(T_{2,2k+1})=
\begin{cases}
-2 & \text{if } r=\tdfrac{1}{4k+1}, \\[1ex]
-2k & \text{if } r=\ndfrac{2k-1}{4k+1},\ \tdfrac{2k}{4k+1}.
\end{cases}
\]
So,
\[
f(r) =
\begin{cases}
-\tdfrac{2}{4k+1} & \text{if } r=\tdfrac{1}{4k+1}, \\
\tdfrac{2(k-2)}{4k+1} & \text{if } r=\ndfrac{2k-1}{4k+1}, \\
\tdfrac{2k}{4k+1} & \text{if } r=\tdfrac{2k}{4k+1}.
\end{cases}
\]
Since $k\ge 3$, $f(r)$ can be negative only when $r=\tdfrac{1}{4k+1}$, and
$f(\tdfrac{1}{4k+1})$ is the minimum.
Thus, $\min_\chi \sigma_1\tau(D_k(K),\chi) =
\min_r\left(2\sigma_{2r}(K)+2f(r)\right) \ge \scrm - \tdfrac{4}{4k+1}$.
This proves (a).
Assuming $\scrm\ge 0$, we have the first statement of (c).

Next, to prove (b) and the second statement of (c),
we compute $f\left(\ttfrac{s}{p}\right)$ when $p=4s\pm 1$.
Let $c = \tnfrac{4k+1}{p}$.
By Proposition~\ref{prop:torus}(a) we have
\[
\sigma_\frac{s}{p}(T_{2,2k+1}) =
-2\left[\frac{s(2k+1)}{p}+\frac{1}{2} \right]
= -2\left[ \frac{cs+1}{2}+\frac{s}{2p}\right]
\]
since $2k+1=\tnfrac{cp+1}{2}$.
Observe that $0<\tdfrac{s}{2p}=\tdfrac{s}{8s\pm 2} \le \tfrac{1}{6}$ 
since $s\ge
1$ and so $8s\pm 2 \ge 6s$.
Since $\tnfrac{cs+1}{2}$ is either an integer or an integer plus
$\tfrac{1}{2}$,
$\sigma_\ttfrac{s}{p}(T_{2,2k+1})=-2\left[\tnfrac{cs+1}{2}\right]$.
Then
\begin{eqnarray*}
f\left(\frac{s}{p}\right) & = &
-2 \left[ \frac{cs+1}{2} \right] + \frac{2s}{p}\left( 1 - \frac{s}{p}
\right) cp \\
& \ge & -(cs +1) + \frac{2cs(p-s)}{p} \\
& = & \frac{cs(p-2s)}{p}-1.
\end{eqnarray*}
If $p=4s+1$, then $\tnfrac{p-2s}{p}=\ndfrac{2s+1}{4s+1}>\frac{1}{2}$. If
$p=4s-1$, then $\ttfrac{s}{p}=\tdfrac{s}{4s-1}>\frac{1}{4}$. Thus we
have
\begin{eqnarray*}
f\left(\frac{s}{p}\right) & > &
\begin{cases}
\dfrac{cs}{2} -1 & \hbox{if } p=4s+1, \\[1ex]
\dfrac{cs}{2} - \dfrac{c}{4} -1 & \hbox{if } p=4s-1.
\end{cases} \\
& \ge & \dfrac{cs}{4} - 1.
\end{eqnarray*}
Note $\tdfrac{cs}{4k+1}=\ttfrac{s}{p}=\tdfrac{s}{4s\pm 1} \ge
\frac{1}{5}$ or
$cs\ge \tnfrac{4k+1}{5}$.  Thus, if $k$ is sufficiently large then
so is $\ttfrac{cs}{4}-1$.
This completes the proof.
\end{proof*}

%%%%%%%%%%%%%%% Subsection %%%%%%%%%%%%%%%%

\subsection{Proofs of Theorem~\ref{thm:rinfinite} and
Theorem~\ref{thm:ribbon}}
\label{sec:proof}
\begin{proof*}[Proof of Theorem~\ref{thm:rinfinite}.]
Let $\scrm = 2\min_{0<r<1} \sigma_r(K)$.
By Lemma~\ref{lem:minmax}(b) there is $k_0>0$ such that
$\sigma_1\tau(D_k(K),\chi_\ttfrac{s}{p}) > |\scrm| +1$
for any $k\ge k_0$.
Let $\scri=\{k\in\bbz \mid k<0 \text{ or } k\ge k_0\}$.  We will show
that $\scri$ is a set satisfying the conclusion of (a).

If $k<0$ then $D_k(K)$ has infinite order in
$\scrg$ by Theorem~\ref{thm:levine}.
Thus, for any $n>0$, $nD_k(K)$ has no metabolizer for the isometric
structure and hence has nonvanishing Gilmer ribbon obstruction.
From now on, assume $k\ge k_0$.

Suppose to the contrary that ${n}D_k(K)$ has vanishing Gilmer ribbon
obstruction for a positive integer $n$.
Then it satisfies the conclusion of Theorem~\ref{thm:gilmer}(b)
for $q=2$.
Let $F$ be the Seifert surface for $D_k(K)$ as depicted in
Figure~\ref{fig:double}.  We take $F_n$ as the boundary connected
sum of $n$ copies of $F$ so that $F_n$ is a Seifert surface of $nD_k(K)$.
Then there is a metabolizer $Z_n$ for the isometric structure
on $H_1(F_n)$ such
that $\tau(nD_k(K),(N^2_n)_p\cap(Z_n\otimes\bbq/\bbz))$ vanishes for all
primes $p$, where $(N^2_n)_p$ is the $p$-primary component of the
kernel, $N^2_n$, of $\varepsilon^2_n\otimes \id_{\bbq/\bbz}$.

Since $\varepsilon^2_n$ is the direct sum of $n$ copies of the map
$\varepsilon^2\!:H_1(F)\to H_1(F)$ corresponding to $D_k(K)$,
$N^2_n$ is the direct sum of $n$ copies of
$N^2=\ker \left(\varepsilon^2\otimes \id_{\bbq/\bbz}\right)$.
Note that $\ker \left(\varepsilon^2\otimes \id_{\bbq/\bbz}\right)\cong
(G^2-(G-I)^2)^{-1}(\bbz\oplus\bbz)
= \left(\begin{smallmatrix}
-1 & -\tdfrac{2k}{4k+1} \\ 2 & 1
\end{smallmatrix}\right)(\bbz\oplus\bbz)$, where
$G$ is the matrix of the isometric structure on $H_1(F)$ with respect to
the basis $\{x,y\}$ as defined in subsection~\ref{sec:rest}.
So $N^2$ is generated by an element
$\left( \tdfrac{1}{4k+1}, 0\right) = x\otimes \tdfrac{1}{4k+1}$.
Thus every character in $(N^2_n)_p$
is a direct sum of characters of the form $x\otimes \ttfrac{s}{p^e}$ and 
$(N^2_n)_p$ is isomorphic to
$(\bbz/p^e)^n$, where $p^e$ is the maximal power of $p$ dividing $4k+1$.
Gilmer~\cite[lemma~2]{gil93} showed that
$|N^2_n|=|N^2_n\cap(Z_n\otimes\bbq/\bbz)|^2$ and hence $|(N^2_n)_p| =
|(N^2_n)_p\cap(Z_n\otimes\bbq/\bbz)|^2$.

Using this and the Gauss-Jordan algorithm,
Livingston and Naik~\cite[proof of Theorem 1.2]{ln99b}
showed that $(N^2_n)_p\cap(Z_n\otimes\bbq/\bbz)$ has an element in
$(\bbz/p^e)^n\cong (N^2_n)_p$ having the
first $n-n_0$ entries equal to $p^{e-1}$
and all the remaining $n_0$ entries divisible by $p^{e-1}$ for some
$n_0\le \ttfrac{n}{2}$. Let $s$ be an integer for which $p=4s\pm 1$.
Multiplying by $s$, we see that $(N^2_n)_p\cap(Z_n\otimes\bbq/\bbz)$ has
an element $\chi$ of the form
$(x\otimes\ttfrac{s}{p},\ldots,x\otimes\ttfrac{s}{p},
x\otimes\ttfrac{s_1}{p},\ldots,x\otimes\ttfrac{s_{n_0}}{p})$,
where $s_i$ can be any integers.
Thus, $\sigma_1\tau(nD_k(K),\chi)  =0$ by the contradiction hypothesis.

On the other hand,
by the additivity of $\sigma_1\tau$, we have
\[
\sigma_1\tau(nD_k(K),\chi)  =
(n-n_0)\sigma_1\tau\left(D_k(K),\chi_\frac{s}{p}\right)
+ \sum_{i=1}^{n_0} \sigma_1\tau\left(D_k(K),\chi_\frac{s_i}{p}\right).
\]
Observe that
$\sigma_1\tau(D_k(K),\chi_\ttfrac{s}{p}) > |\scrm| +1 >
|\scrm|+ \tdfrac{4}{4k+1}$
for any $k\ge k_0$.  Now by Lemma~\ref{lem:minmax}(a),
\[
\sigma_1\tau(nD_k(K),\chi) >
(n-n_0) \left( |\scrm| + \frac{4}{4k+1} \right)
+ n_0 \left( \scrm - \frac{4}{4k+1} \right) \ge 0
\]
since $n-n_0\ge n_0$.  So $\sigma_1\tau(nD_k(K),\chi) >0$. This is a
contradiction, proving (a).

Next, assume $\scrm\ge 0$ and $k\ge 3$.  Let $\chi$ denote the
character as given above again.
If $4k+1$ is a composite number, i.e.~
$\tnfrac{4k+1}{p}>1$, then $\sigma_1\tau(D_k(K),\chi) > 0$ by
Lemma~\ref{lem:minmax}(c) since none of
$\ttfrac{s}{p}$ and $\ttfrac{s_i}{p}$ is $\tdfrac{1}{4k+1}$.
Now assume $4k+1=p$. Then $c=1$ and $s=k$, where $c$ and $s$ are those in
Lemma~\ref{lem:minmax}(c), and
$\sigma_1\tau(D_k(K),\chi_\ttfrac{k}{p})>k-2>\tdfrac{4}{4k+1}$.  We
now have
\[
\sigma_1\tau(nD_k(K),\chi)>(n-n_0)\frac{4}{4k+1} + n_0\frac{-4}{4k+1}
\ge 0.
\]
This proves the first part of (b).

For $k=1,2$, an elementary computation shows:
\[
\begin{array}{|c||c|c|c|c|c|c|}
\hline
k &
\multicolumn{2}{c|}{1} & \multicolumn{4}{c|}{2} \\
\hline
r\rule{0em}{1.1em} &
\frac{1}{5} & \frac{2}{5} &
\frac{1}{9} & \frac{2}{9} & \frac{3}{9} & \frac{4}{9} \\[.5ex]
\hline
\frac{1}{2}\sigma_1\tau(D_k(K),\chi_r) - \sigma_{2r}(K)\rule{0em}{1.1em}
& -\frac{2}{5} & \frac{2}{5} &
-\frac{2}{9} & \frac{10}{9} & 0 & \frac{4}{9} \\[.5ex]
\hline
\end{array}
\]
So if $\sigma_{2r}(K) > 0$ for $r=\frac{1}{5}, \frac{1}{9},
\frac{3}{9}$, then $\sigma_1\tau(D_k(K),\chi_\tdfrac{s}{4k+1})$
are all positive for $k=1,2$.  Note that
$\sigma_\frac{6}{9}(K)=\sigma_{1-\frac{1}{3}}(K)
=\sigma_\frac{1}{3}(K)$.
This completes the proof.
\end{proof*}

\begin{proof*}[Proof of Theorem~\ref{thm:ribbon}.]
The exact same proof of Theorem~\ref{thm:independent} works here
by applying all the counter-parts for the ribbon
case: For instance, Corollary~\ref{cor:ralexander} instead of
Theorem~\ref{thm:alexander}.
\end{proof*}

%%%%%%%%%%%%%%%%%% Acknowledgments %%%%%%%%%%%%%%%%%%%%%%%%%

\begin{acknowledgements}
The author thanks Chuck Livingston for his suggestions and
discussions.
He also thanks the referee of the first version of the paper for
many useful comments.
\end{acknowledgements}

%%%%%%%%%%%%%%%%%% References %%%%%%%%%%%%%%%%%%%%%%%%%

\end{document}